\documentclass[final,1p,times]{elsarticle}

\usepackage{amssymb}
 \usepackage{amsthm}
\usepackage{amscd}
\usepackage{amsmath}
\usepackage{amsfonts}
\usepackage{amssymb}
\usepackage{graphicx}
\usepackage{mathrsfs}
\usepackage{titletoc}

\numberwithin{equation}{section}

\newcommand{\ra}{\rightarrow}
\newcommand{\p}{\partial}
\newcommand{\f}{\frac}

\newcommand{\be}{\begin{equation}}
\renewcommand{\ra}{\rightarrow}
\newcommand{\ee}{\end{equation}}
\newcommand{\bea}{\begin{eqnarray}}
\newcommand{\eea}{\end{eqnarray}}
\newcommand{\bna}{\begin{eqnarray*}}
\newcommand{\ena}{\end{eqnarray*}}

\renewcommand{\le}{\left}
\newcommand{\ri}{\right}

\journal{***}

\begin{document}

\begin{frontmatter}

\title{Quantization for an elliptic equation with
critical exponential growth on compact Riemannian surface without
boundary}

\author{Yunyan Yang}
 \ead{yunyanyang@ruc.edu.cn}
\address{ Department of Mathematics,
Renmin University of China, Beijing 100872, P. R. China}

\begin{abstract}
  In this paper, using blow-up analysis, we prove a quantization result for an elliptic equation with critical exponential growth
  on compact Riemannian surface without boundary. Similar results for Euclidean space were obtained by
  Adimurthi-Struwe \cite{Adi-Stru}, Druet \cite{Druet}, Lamm-Robert-Struwe \cite{L-R-S}, Martinazzi \cite{Mart},
  Martinazzi-Struwe \cite{Mar-Stru}, and Struwe \cite{Struwe} respectively.
\end{abstract}


\end{frontmatter}

\titlecontents{section}[0mm]
                       {\vspace{.2\baselineskip}}
                       {\thecontentslabel~\hspace{.5em}}
                        {}
                        {\dotfill\contentspage[{\makebox[0pt][r]{\thecontentspage}}]}
\titlecontents{subsection}[0mm]
                       {\vspace{.2\baselineskip}}
                       {\thecontentslabel~\hspace{.5em}}
                        {}
                       {\dotfill\contentspage[{\makebox[0pt][r]{\thecontentspage}}]}

\setcounter{tocdepth}{1}

\setcounter{tocdepth}{1}

\tableofcontents

\section{Introduction and main results}
 Let $(\Sigma,g)$ be a compact Riemannian surface without boundary,
 $W^{1,2}(\Sigma,\mathbb{R})$ be the usual Sobolev space, namely the completion
 of $C^\infty(\Sigma,\mathbb{R})$ under the norm
 $$\|u\|_{W^{1,2}(\Sigma,\mathbb{R})}=\le(\int_{\Sigma}\le(|\nabla_gu|^2+u^2\ri)dv_g\ri)^{1/2},$$
 where $\nabla_gu$ denotes the gradient of $u$ and $dv_g$ denotes the volume element with respect to the Riemannian metric $g$.
  Let $f_k:\Sigma\times[0,\infty)\ra
 \mathbb{R}$ be a sequence of functions satisfying the
 following hypotheses:\\

 \noindent(H1)\,\,$f_k(x,0)=0$, and $f_k(x,t)>0$ for all
 $k$, all $x\in \Sigma$, and all $t>0$;\\
 (H2)\,\,$f_k\in C^2(\Sigma\times[0,+\infty))$ for each $k$ and $f_k\ra f_\infty$ in
$C^1_{\rm
 loc}(\Sigma\times[0,+\infty))$ as $k\ra\infty$;\\
 (H3)\,\,for any $\nu>0$, there exists a constant
 $C_\nu>0$ such that for all $k$, all $x\in\Sigma$, and all $t>0$,
 $$F_k(x,t)\leq \nu tf_k(x,t)+C_\nu,$$
 where
 $$F_k(x,t)=\int_0^tf_k(x,s)ds$$
 is the primitive of $f_k(x,t)$;\\
 (H4)\,\,$f_k^\prime(x,t)/(tf_k(x,t))\ra 2$ as $t\ra+\infty$
 uniformly in $k\in\mathbb{N}$ and in $x\in\Sigma$, where $f_k^\prime$ is the
 derivative of $f_k$ with respect to $t$, moreover there exists a constant $C$ such that
 $|\nabla_gf_k(x,t)|\leq C(1+f_k(x,t))$ for all
 $(x,t)\in\Sigma\times\mathbb{R}$;\\
 (H5)\,\,there exist $\psi$, a continuous function with $\psi(0)=0$, $t_0>0$,
 and $k_0>0$, such that
  $$|f_k(x,t)/f_k(y,t)-1|\leq\psi({d}_g(x,y))$$
  for all $t\geq t_0$, all $k>k_0$, and all $x,y\in\Sigma$, where $d_g(\cdot,\cdot)$ denotes the
  geodesic distance between two points of $\Sigma$.\\

  By (H4) we have $f_k(x,t)=f_k(x,t_0)e^{(1+o(1))(t^2-t_0^2)}$ for any given $t_0>0$, where $o(1)\ra 0$
  as $t\ra\infty$ uniformly in $x\in\Sigma$.
  In view of the Trudinger-Moser embedding
 \cite{Fontana,Moser,Pohozaev,Trudinger}, we say that $f_k(x,t)$ is of critical exponential
 growth with respect to $t$. A typical example satisfying (H1)-(H5) is
 \be\label{sequence}f_k(x,t)=\lambda_k te^{t^2},\ee
 where $\lambda_k$ is a sequence of positive real numbers such that
 $\lambda_k\ra\lambda_\infty$
 as $k\ra\infty$.
 Suppose that for each $k\in\mathbb{N}$ we have a smooth function $u_k\geq
 0$ satisfying the equation
 \be\label{e1}\Delta_gu_k+\tau_ku_k=f_k(x,u_k)\quad{\rm in}\quad\Sigma,\ee
 where $\Delta_g$ is the Laplace-Beltrami operator, $\tau_k$ is a sequence of smooth
 functions such that
 \be\label{tauk}\tau_k\ra\tau_\infty\,\,\,{\rm in}\,\,\,
 C^0(\Sigma,\mathbb{R}),\,\,\,
 \tau_\infty(x)>0\,\,\,{\rm for\,\,\,all}\,\,\, x\in \Sigma.\ee
 Clearly $u_k$ is a critical point of the functional
 \be\label{functional}J_k(u)=\f{1}{2}\int_\Sigma\le(|\nabla_gu|^2+\tau_ku^2\ri)dv_g-\int_\Sigma F_k(x,u)dv_g\ee
 on the Sobolev space $W^{1,2}(\Sigma,\mathbb{R})$. The existence of
 nonnegative solutions to equation (\ref{e1}) in case that $\tau_k$ is a positive real number
 was studied by  Zhao and the author
 \cite{Yang-Zhao} by using variational methods. More explicitly, assuming that
 $\lambda_\tau=\lambda_\tau(\Sigma)$ is the first eigenvalue of
 the operator $\Delta_g+\tau$, where $\tau>0$ is a constant, we proved that
 the equation $\Delta_gu+\tau u=\lambda ue^{u^2}$ has a nonnegative solution if
 $\lambda<\lambda_{\tau}$. The aim of this
 paper is to study the quantization problem for equation (\ref{e1}).
 Precisely we shall prove the following result.\\

\noindent{\bf Theorem 1.1} {\it Let $(\Sigma,g)$ be a compact
Riemannian surface without boundary. Suppose that $u_k\geq 0$ is a
sequence of smooth solutions to equation (\ref{e1}),
 where $\tau_k$ is a sequence of smooth functions satisfying
 (\ref{tauk}),
 and $f_k$ is a sequence of functions
 satisfying {\rm (H1)-(H5)}. Let $J_k$ be as in (\ref{functional}).
 If $J_k(u_k)\ra \beta$ as $k\ra\infty$ for some
 $\beta\in\mathbb{R}$, then there exists a nonnegative solution $u_\infty\in
 C^1(\Sigma,\mathbb{R})$ of the equation
 \be\label{3'}\Delta_gu_\infty+\tau_\infty u_\infty=f_\infty(x,u_\infty(x))\quad{\rm in}\quad
 \Sigma,\ee
 and there exists $N\in \mathbb{N}$ such that
 $J_k(u_k)=J_\infty(u_\infty)+2\pi N+o(1)$, where $o(1)\ra 0$ as $k\ra\infty$. Here $J_\infty$ is also as in
 (\ref{functional}), where $\tau_k$, $F_k$ are replaced by $\tau_\infty$
 and $F_\infty$ respectively.
 If $N=0$, $u_k\ra u_\infty$ strongly in $W^{1,2}(\Sigma,\mathbb{R})$ and in
 fact in $C^1(\Sigma,\mathbb{R})$.}\\

 Several works were devoted to prove analogues of Theorem 1.1. In
 \cite{Adi-Stru}, Adimurthi and Struwe considered a sequence of
 solutions $u_k$ to the equation
 \be\label{eucl}\le\{\begin{array}{lll}-\Delta_{\mathbb{R}^2} u_k=f_k(x,u_k)\,\,\,{\rm in}\,\,\,
 \Omega\subset\mathbb{R}^2\\[1.5ex]u_k>0\,\,\,{\rm in}\,\,\,\Omega,\,\,\,u_k=0\,\,\,{\rm on}\,\,\,\p\Omega,\end{array}\ri.\ee
 where $f_k(x,t)=te^{\varphi_k(t)}$, $0\leq\varphi_k^{\prime\prime}(t)\leq
 2$ for $t\geq t_0$ and $\varphi_k^\prime(t)/t\ra 2$ as $t\ra\infty$
 uniformly in $k$. Such a sequence of functions $f_k$ satisfies (H1)-(H5) in case that the Riemannian surface $(\Sigma,g)$
 is replaced by a smooth bounded domain of $\mathbb{R}^2$. Assuming that
 $$J_k(u_k)=\f{1}{2}\int_\Omega |\nabla_{\mathbb{R}^2} u_k|^2dx-
 \int_\Omega F_k(x,u_k)dx\ra \beta$$ for $0\leq\beta<4\pi$ and that
 the limit equation does not admit any positive solution with energy
 less than $2\pi$, they proved that either $u_k\ra u_\infty$ strongly in $W_0^{1,2}(\Omega)$ and
 $u_\infty$ has energy $\beta$, or $u_k\rightharpoonup 0$
 weakly in $W_0^{1,2}(\Omega)$ and $u_k$ develops one blow-up
 point carrying the energy $2\pi$. This quantization result was surprisingly refined by Druet
 \cite{Druet} to the case of all $\beta\in\mathbb{R}$ and general
 nonlinearities of uniform critical growth, analogous to that of the current paper.
 (Blow-up analysis for equation (\ref{eucl}) with similar nonlinearity was
 also considered by Adimuthi and Druet \cite{A-D}.) The key point in \cite{Druet} is
 the gradient estimate (\cite{Druet}, Proposition 2), through which Druet studied the energy of ${\varphi}_k$,
 the spherical average of $u_k$ with respect to blow-up
 points,
 instead of $u_k$ itself. Thus he transformed the quantization problem for $u_k$ to the quantization problem
 for ${\varphi}_k$, which depends only on analysis on certain ordinary differential
 equation and is comparatively easy to be handled.
 Shortly after, using similar idea, Struwe
 \cite{Struwe} succeeded to get a quantization result for a forth order elliptic
 equation
 $$\le\{\begin{array}{lll}-\Delta_{\mathbb{R}^4}^2 u_k=\lambda_ku_ke^{2u_k^2}\,\,\,{\rm in}\,\,\,
 \Omega\subset\mathbb{R}^4\\[1.5ex]u_k>0\,\,\,{\rm in}\,\,\,\Omega,\,\,\,u_k=\Delta_{\mathbb{R}^4}u_k=0
 \,\,\,{\rm on}\,\,\,\p\Omega,\end{array}\ri.$$
 where $0<\lambda_k\ra 0$ as $k\ra\infty$, and $u_k\rightharpoonup
 0$ weakly in $W^{2,2}(\Omega)$. Also Lamm, Robert and Struwe \cite{L-R-S} proved
 a quantization result for the evolution of equation
 (\ref{eucl}), where $f_k$ is as in (\ref{sequence}). A recent inspiring work of
 Martinazzi and Struwe \cite{Mar-Stru} states the following: Let
 $\Omega\subset\mathbb{R}^{2m}$ be a smooth bounded domain, $u_k$
 be a sequence of positive solutions to the equation $-\Delta_{\mathbb{R}^2}u_k=\lambda_k u_k e^{mu_k^2}$
 subject to Dirichlet boundary conditions, where $0<\lambda_k\ra 0$
 and $u_k\rightharpoonup 0$ weakly in $W^{m,2}(\Omega)$. Assuming $\Lambda=\lim_{k\ra\infty}\int_\Omega
 u_k(-\Delta_{\mathbb{R}^{2m}})^mu_kdx<\infty$, they proved that
 $\Lambda$ is an integer multiple of $\Lambda_1=(2m-1)!{\rm
 vol}(\mathbb{S}^{2m})$, the total $Q$-curvature of the standard
 $2m$-dimensional sphere. In view of the Trudinger-Moser embedding
 for the space $W_0^{1,n}(\Omega)$, where $n\geq 3$ and $\Omega\subset\mathbb{R}^n$
 is a smooth bounded domain, one may ask how about the equation
 \be\label{nD}\le\{\begin{array}{lll}-\Delta_n u_k=\lambda_ku_k^{\f{1}{n-1}}e^{u_k^{\f{n}{n-1}}}\,\,\,{\rm in}\,\,\,
 \Omega\\[1.5ex]u_k\geq 0\,\,\,{\rm in}\,\,\,\Omega,\,\,\,u_k=0
 \,\,\,{\rm on}\,\,\,\p\Omega.\end{array}\ri.\ee
 Up to now only an energy inequality has been obtained by Adimurthi and the
 author \cite{Adi-Yang}. Concerning the quantization for
 equation (\ref{nD}), we have a long way to go. For other works related to
 this kind of
 quantization problems we refer the reader to
 \cite{Mar-Stru,Struwe} and the references therein.\\

 For the proof of Theorem 1.1, we follow the lines of \cite{Druet,L-R-S,Mar-Stru,Struwe}.
 Firstly we use a pointwise estimate on $u_k$ to find all separate blow-up
 points. Specifically we need to deal  carefully with the term $\tau_k u_k$, which does not appear in the
 Euclidean case. Secondly we establish a gradient estimate for $u_k$.
 This permits us to compare $u_k$ with its spherical
 average with respect to blow-up points. Finally
 we get the quantization result, where we should deal with the extra term $\tau_k u_k$ again.
 For calculations near blow-up points we prefer to choose
 isothermal coordinates instead of normal coordinates. The advantage of such coordinates is
 that both the Laplace-Beltrami operator $\Delta_g$ and the gradient
 operator $\nabla_g$ have simple expressions.\\

 The remaining part of this paper is organized as follows. In the
 next section we prove a simple property of the weak convergence of $u_k$.
 In Section 3, we locate the blow-up points of $u_k$ and
 describe the asymptotic behavior of $u_k$ near those points. In
 Section 4 we derive a gradient estimate on $u_k$.
 We shall prove quantization results for $u_k$ near the blow-up points in Section
 5, and complete the proof of Theorem 1.1
 in Section 6.\\

 Throughout this paper we often denote various constants independent of $k$ by the same
 $C$. In addition, we do not distinguish between sequence and subsequence or points and sequence sometimes.
 The
 reader can easily recognize it from the context.

 \section{Weak convergence}

 In this section, we let $u_k\geq 0$ be a sequence of solutions to
 equation (\ref{e1}) verifying that
 \be\label{beata}J_k(u_k)\ra\beta\,\,\, {\rm as}\,\,\, k\ra\infty\,\,\,
 {\rm for\,\, some}\,\,\, \beta\in\mathbb{R},\ee
 where $J_k$ is defined in (\ref{functional}).
 Testing equation (\ref{e1}) by
$u_k$, we have
\be\label{0.7}
\int_\Sigma\le(|\nabla_gu_k|^2+\tau_ku_k^2\ri)dv_g=\int_\Sigma
u_kf_k(x,u_k)dv_g. \ee It follows from (\ref{beata}) that
$$\int_\Sigma\le(|\nabla_gu_k|^2+\tau_ku_k^2\ri)dv_g=2\beta+2\int_\Sigma F_k(x,u_k)dv_g+o(1).$$
Hence
$$\int_\Sigma u_kf_k(x,u_k)dv_g=2\beta+2\int_\Sigma F_k(x,u_k)dv_g+o(1).$$
 If $f_k$ satisfies the hypotheses (H1)-(H4), then we have
\be\label{bdd} \int_\Sigma u_kf_k(x,u_k)dv_g\leq C\ee for some
constant $C$. In view of (\ref{tauk}), it follows from (\ref{0.7})
and (\ref{bdd}) that $u_k$ is bounded in
$W^{1,2}(\Sigma,\mathbb{R})$. Hence there exists some $u_\infty\in
W^{1,2}(\Sigma,\mathbb{R})$ such that up to a subsequence,
$u_k\rightharpoonup u_\infty$ weakly in
$W^{1,2}(\Sigma,\mathbb{R})$, $u_k\ra u_\infty$ strongly in
$L^2(\Sigma,\mathbb{R})$, and $u_k\ra u_\infty$ a.e. in $\Sigma$.
Similarly to \cite{Druet}, we then get that
\be\label{F}\lim_{k\ra\infty}\int_\Sigma F_k(x,u_k)dv_g=\int_\Sigma
F_\infty(x,u_\infty)dv_g\ee
 that  $u_\infty$ is a weak
solution of (\ref{3'}), and that  $u_\infty\in
C^1(\Sigma,\mathbb{R})$. In conclusion we obtained an analogue of (\cite{Druet}, Lemma1),
namely\\

\noindent{\bf Lemma 2.1} {\it Let $f_k$ be a sequence of functions
satisfying (H1)-(H4). Let $u_k\geq 0$ be a sequence of
solutions to (\ref{e1}), where $\tau_k$ is
as defined in (\ref{tauk}). If (\ref{beata}) holds, then $u_k$ is
bounded in $W^{1,2}(\Sigma,\mathbb{R})$, and thus, up to a subsequence,
$u_k\rightharpoonup u_\infty$ weakly in
$W^{1,2}(\Sigma,\mathbb{R})$, where $u_\infty\in
C^1(\Sigma,\mathbb{R})$ is a solution to (\ref{3'}).
Also, there holds
\be\label{0.6}\lim_{k\ra\infty}\int_\Sigma\le(|\nabla_gu_k|^2+\tau_ku_k^2\ri)dv_g=
2\beta+2\int_\Sigma F_\infty(x,u_\infty)dv_g.\ee }

 \section{Multibubble analysis}

 In this section we shall use point wise estimate to find blow-up
 points of a sequence of solutions to the equation
 (\ref{e1}). This technique was first used  by Druet \cite{Druet} to deal with blow-up analysis for
 solutions to the equation (\ref{eucl}). Assume $u_k\geq 0$ is a sequence of solutions to the equation (\ref{e1})
 and (\ref{beata}) holds. From (\ref{0.7}) and (\ref{bdd}) we can
 find some constant $C$ such that
\be\label{energy-bounded}\int_\Sigma
\le(|\nabla_gu_k|^2+\tau_ku_k^2\ri)dv_g\leq C.\ee
 Then the Sobolev embedding theorem implies that for any $p>1$ there is
 some constant $C$ such that
 \be\label{Lp-bdd}\int_\Sigma u_k^pdv_g\leq C.\ee
 These two properties are very important during the process of exhausting blow-up points.
 Precisely we have the following proposition which is
 analogous to (\cite{Druet}, Proposition 1), (\cite{L-R-S}, Theorem 4.2), (\cite{Mart}, Theorem 1 in the case
 $m=1$) and (\cite{Adi-Yang}, Proposition 3.1).
  \\

 \noindent{\bf Proposition 3.1} {\it Let $(\Sigma,g)$ be a compact Riemannian surface without boundary,
 $(f_k)$ be a sequence of functions satisfying the hypotheses {\rm (H1)-(H5)}, and $(u_k)$ be a sequence of smooth
 nonnegative solutions to (\ref{e1}) such that (\ref{beata}) holds.
 Assume that $\max_\Sigma u_k\ra+\infty$
 as $k\ra\infty$. Then there exists $N\in\mathbb{N}\setminus\{0\}$, and up to a subsequence, there exist
 $N$ sequences of points $x_{i,k}\ra x_i^\ast\in\Sigma$ and of positive real numbers $r_{i,k}\ra 0$ as $k\ra\infty$, where
 $r_{i,k}$ is defined by
 \be\label{scal}r_{i,k}^{-2}=u_k(x_{i,k})f_k(x_{i,k},u_k(x_{i,k})),\ee
 such that the following hold:\\

 \noindent$(i)$ For any $i=1,2,\cdots,N$, take an isothermal coordinate system
 $(U_i,\phi_i;\{x^1,x^2\})$ near $x_i^\ast$, where $U_i\subset\Sigma$ is a
 neighborhood of $x_i^\ast$, $\phi_i: U_i\ra \Omega_i\subset\mathbb{R}^2$ is a diffeomorphism
 and $\phi_i(x_i^\ast)=(0,0)$. If we define
 \be\label{blow-up-funct}\eta_{i,k}(x)=u_k(x_{i,k})(\widetilde{u}_k(\widetilde{x}_{i,k}+r_{i,k}x)-u_k(x_{i,k}))\ee
 for all $x\in\Omega_{i,k}=\{x\in\mathbb{R}^2: \widetilde{x}_{i,k}+r_{i,k}\in\Omega_i\}$,
 where
 $\widetilde{x}_{i,k}=\phi_i(x_{i,k})$ and $\widetilde{u}_k=u_k\circ\phi_i^{-1}$,
 then there holds
 $$\eta_{i,k}(x)\ra\eta_\infty(x)=\log\f{1}{1+|x|^2/4}\quad{\rm in}\quad C^1_{\rm loc}(\mathbb{R}^2);$$
 $(ii)$ For any $1\leq i\not=j\leq N$, there holds
 $$\f{d_g(x_{i,k},x_{j,k})}{r_{i,k}}\ra+\infty,
 \quad{\rm as}\quad k\ra\infty,$$ where $d_g(\cdot,\cdot)$
 denotes the geodesic distance between two points of $\Sigma$;\\
 $(iii)$ Define $R_{N,k}(x)=\min_{1\leq i\leq N}d_g(x,x_{i,k})$ for $x\in\Sigma$, then there exists a constant $C>0$ such that
 $$R^2_{N,k}(x)u_k(x)f_k(x,u_k(x))\leq C$$
 uniformly in $x\in \Sigma$ and $k\in \mathbb{N}$.\\

\noindent Moreover, given any sequence of points $(x_{N+1,k})$, it
is impossible to extract a new subsequence from the previous
 one such that $(i)-(iii)$ hold with the sequences $(x_{i,k})$,
 $i=1,\cdots,N+1$.\\

Finally, we have $u_k\ra u_\infty$ in $C^1_{\rm
loc}(\Sigma\setminus\mathcal{S})$ as $k\ra\infty$, where
$\mathcal{S}=\{x_1^\ast,\cdots,x_N^\ast\}$, and
 $u_\infty$ is given in Lemma 2.1.}\\

 {\it Proof.} Similarly to \cite{Druet,L-R-S,Mart,Adi-Yang}, we prove the proposition
  by several steps as follows.\\

 {\it Step $1$. The first bubble.}\\

 Assume $u_k(x_k)=\max_\Sigma u_k$. If $u_k(x_k)$ is
 bounded, applying elliptic estimates to equation (\ref{e1}), we then have
 $u_k\ra u_\infty$ in $C^1(\Sigma,\mathbb{R})$, where
 $u_\infty$ is given by Lemma 2.1. Hereafter we assume
 $u_k(x_k)\ra+\infty$. Set
 \be\label{rk}r_k^{-2}=u_k(x_k)f_k(x_k,u_k(x_k)).\ee
 It is clear that $r_k\ra 0$ as $k\ra\infty$.

 Assume $x_k\ra x^\ast$ as $k\ra\infty$. Take an isothermal coordinate system $(U,\phi;\{x^1,x^2\})$
 near $x^\ast$, where $U\subset\Sigma$ is a
 neighborhood of $x^\ast$, $\phi: U\ra \Omega\subset\mathbb{R}^2$ is a diffeomorphism
 and $\phi(x^\ast)=(0,0)$.
 In such a coordinate system, the metric $g$ can be represented by
 $$g=e^{\psi}(d{x^1}^2+d{x^2}^2)$$
 for some smooth function $\psi:\Omega\ra\mathbb{R}$ with
 $\psi(0,0)=0$. It follows that
 \be\label{laplace}\nabla_g=e^{-\psi}\nabla_{\mathbb{R}^2},\quad \Delta_g=-e^{-\psi}
 \Delta_{\mathbb{R}^2},\ee
 where $\nabla_{\mathbb{R}^2}$ and $\Delta_{\mathbb{R}^2}$ denote the usual gradient operator
 and the Laplace operator of $\mathbb{R}^2$ respectively. The
 existence of isothermal coordinate system on Riemannian surface is a well-known fact in Riemannian geometry,
 see for example \cite{Wu}.
 Define
 \be\label{vk}v_k(x)=\f{\widetilde{u}_k(\widetilde{x}_k+r_kx)}{u_k(x_k)}\ee
 for $x\in\Omega_k=\{x\in\mathbb{R}^2:\widetilde{x}_k+r_kx\in\Omega\}$,
 where $\widetilde{u}_k=u_k\circ\phi^{-1}$,
 $\widetilde{x}_k=\phi(x_k)$.
 It follows from (\ref{e1}), (\ref{rk}) and (\ref{laplace}) that $v_k$ satisfies the following equation
 \be\label{11}
 -\Delta_{\mathbb{R}^2}v_k(x)=e^{\psi(\widetilde{x}_k+r_kx)}\f{\widetilde{f}_k(\widetilde{x}_k+r_kx,
 \widetilde{u}_k(\widetilde{x}_k+r_kx))}{u_k^2(x_k)f_k(x_k,u_k(x_k))}
 -e^{\psi(\widetilde{x}_k+r_kx)}r_k^2\widetilde{\tau}_k(\widetilde{x}_k+r_kx) v_k(x)\ee
 on $\Omega_k$,
 where
 $\widetilde{f}_k(\widetilde{x}_k+r_kx,t)={f}_k(\phi^{-1}(\widetilde{x}_k+r_kx),t)$.
 Note that $u_k(x_k)=\max_\Sigma u_k$ and $\Omega_k\ra\mathbb{R}^2$ as $k\ra\infty$. It follows from (\ref{vk})
 that $v_k$ is uniformly bounded in $\mathbb{B}_R(0)$ for any fixed
 $R>0$. Since $\psi$ is smooth, $\psi(0,0)=0$, $\widetilde{x}_k\ra(0,0)$
 and $r_k\ra 0$ as $k\ra\infty$, $e^{\psi(\widetilde{x}_k+r_kx)}$
 is also uniformly bounded in $\mathbb{B}_R(0)$ for any fixed
 $R>0$. Furthermore $e^{\psi(\widetilde{x}_k+r_kx)}\ra 1$ locally uniformly in $\mathbb{R}^2$
 as $k\ra\infty$. By (H4) and (H5), we have for all
 $x\in \Omega_k$ and all $k$
 \be\label{f-bdd}\f{\widetilde{f}_k(\widetilde{x}_k+r_kx,
 \widetilde{u}_k(\widetilde{x}_k+r_kx))}{f_k(x_k,u_k(x_k))}\leq
 C.\ee
 All these estimates together with (\ref{tauk}) lead to
 $$\|-\Delta_{\mathbb{R}^2}v_k\|_{L^\infty(\mathbb{B}_R(0))}\ra 0\,\,\,{\rm
 as}\,\,\,
 k\ra\infty,\,\,\forall R>0.$$
  Applying elliptic estimates to (\ref{11}), one gets $v_k\ra
 v_\infty$ in $C^1_{\rm loc}(\mathbb{R}^2)$, where $v_\infty$
 satisfies
 $$\le\{\begin{array}{lll}
 -\Delta_{\mathbb{R}^2}v_\infty=0\,\,\,{\rm in}\,\,\,\mathbb{R}^2\\[1.5ex]
 v_\infty(0)=1=\max_{\mathbb{R}^2}v_\infty.
 \end{array}\ri.$$
 The Liouville theorem for harmonic functions then leads to $v_\infty\equiv
 1$. Therefore
 \be\label{C-V}v_k\ra 1\quad{\rm in}\quad C^1_{\rm loc}(\mathbb{R}^2).\ee
 Now we set
 $$\eta_k(x)=u_k(x_k)(\widetilde{u}_k(\widetilde{x}_k+r_kx)-u_k(x_k)).$$
 In view of (\ref{e1}), $\eta_k$ satisfies
 \bea{\nonumber}-\Delta_{\mathbb{R}^2}\eta_k(x)&=&e^{\psi(\widetilde{x}_k+r_kx)}\f{\widetilde{f}_k(\widetilde{x}_k+r_kx,
 \widetilde{u}_k(\widetilde{x}_k+r_kx))}{f_k(x_k,u_k(x_k))}\\[1.5ex]&&
 -e^{\psi(\widetilde{x}_k+r_kx)}\widetilde{\tau}_k(\widetilde{x}_k+r_kx) r_k^2u_k^2(x_k)v_k(x),\,\,\, x\in
 \Omega_k.\label{eta}\eea
 We claim that
 \be\label{rk-0}r_ku_k^p(x_k)\ra 0 \quad{\rm as}\quad k\ra\infty,\quad \forall
 p>1.\ee Actually,
 it is clear that there exists some constant $c>0$ depending only on
 the diffeomorphism $\phi$ such that for any fixed $R>0$ and all
 large $k$
 \be\label{ball}B_{c^{-1}Rr_k}(x_k)\subset \phi^{-1}\le(\mathbb{B}_{Rr_k}(\widetilde{x}_k)\ri)\subset
 B_{cRr_k}(x_k).\ee
 Here and throughout this paper we denote the geodesic ball centered at $x\in\Sigma$ with radius $r$ by
 $B_r(x)$, while the Euclidean ball centered at $x\in\mathbb{R}^2$ with radius $r$ by $\mathbb{B}_r(x)$.
 This together with (\ref{C-V}), the mean value theorem for integral and the H\"older
 inequality leads to
 \bea{\nonumber}
 r_ku_k^p(x_k)&=&\f{r_k}{\pi}\int_{\mathbb{B}_1(0)}u_k^p(x_k)dx\\{\nonumber}
 &=&(1+o(1))\f{r_k}{\pi}\int_{\mathbb{B}_1(0)}\widetilde{u}_k^p(\widetilde{x}_k+r_kx)dx\\{\nonumber}
 &\leq&(1+o(1))\f{r_k}{\pi^{1/3}}\le(\int_{\mathbb{B}_1(0)}\widetilde{u}_k^{3p}(\widetilde{x}_k+r_kx)dx\ri)^{1/3}\\
 &\leq&(1+o(1))\f{r_k^{1/3}}{\pi^{1/3}}\le(\int_{{B}_{cr_k}(x_k)}{u}_k^{3p}dv_g\ri)^{1/3},{\label{vrg}}
 \eea
 where $o(1)\ra 0$ as $k\ra\infty$ for any fixed $p>1$. In view of (\ref{Lp-bdd}),
 our claim (\ref{rk-0}) follows from (\ref{vrg}) immediately.

 For any fixed $R>0$ we let $\eta_k^{(1)}$ be a solution to the equation
 \be\label{312}\le\{\begin{array}{lll}
 -\Delta_{\mathbb{R}^2}\eta_k^{(1)}=-\Delta_{\mathbb{R}^2}\eta_k\,\,\,{\rm in}\,\,\,\mathbb{B}_R(0)\\[1.5ex]
 \eta_k^{(1)}=0\quad {\rm on}\quad{\p\mathbb{B}_R(0)}.
 \end{array}\ri.\ee
 In view of (\ref{eta}), we have by (\ref{f-bdd}) and (\ref{rk-0})
 that $\Delta_{\mathbb{R}^2}\eta_k$ is bounded in $L^\infty_{\rm
 loc}(\mathbb{R}^2)$. Applying elliptic estimates to (\ref{312}), we have
 \be\label{bou}\eta_k^{(1)}\ra \eta_\infty^{(1)}\quad{\rm in}\quad
 C^1(\mathbb{B}_{R}(0)).\ee
  Let $\eta_k^{(2)}=\eta_k-\eta_k^{(1)}$. Then $\eta_k^{(2)}$
  satisfies
 \be\label{hm}-\Delta_{\mathbb{R}^2}\eta_k^{(2)}=0\quad{\rm in}\quad
 \mathbb{B}_R(0).\ee
 It follows from (\ref{bou}) and $\eta_k\leq 0$ that there exists some constant $C$  such that
 $\eta_k^{(2)}(x)\leq C$ for all $k$ and all $x\in \mathbb{B}_R(0)$.
 Applying the Harnack inequality to (\ref{hm}), we conclude that
 $\eta_k^{(2)}$ is uniformly bounded on $\mathbb{B}_{R/2}(0)$. Hence
 $\eta_k$ is also uniformly bounded in $\mathbb{B}_{R/2}(0)$. Applying
 elliptic estimates to (\ref{eta}), we obtain
 $$\eta_k\ra \eta_\infty\quad {\rm in}\quad
 C^1(\mathbb{B}_{R/4}(0)).$$
 This together with (H4), (H5) and (\ref{C-V}) gives
 \be\label{b-c}\f{\widetilde{f}_k(\widetilde{x}_k+r_kx,
 \widetilde{u}_k(\widetilde{x}_k+r_kx))}{f_k(x_k,u_k(x_k))}=(1+o(1))e^{(2+o(1))\eta_\infty}
 \ee
 for all $x\in\mathbb{B}_{R/4}(0)$,
 where $o(1)\ra 0$ as $k\ra\infty$ uniformly in $x\in
 \mathbb{B}_{R/4}(0)$. Inserting (\ref{rk-0}) and (\ref{b-c}) into
 (\ref{eta}) and noting that $R>0$ is arbitrary  we obtain
 \be\label{bubble}\le\{\begin{array}{lll}
 -\Delta_{\mathbb{R}^2}\eta_\infty=e^{2\eta_\infty}\,\,\,{\rm in}\,\,\,\mathbb{R}^2\\[1.5ex]
 \eta_\infty(0)=0=\max_{\mathbb{R}^2}\eta_\infty.
 \end{array}\ri.\ee
 Moreover, using (\ref{bdd}), (\ref{rk}), (\ref{C-V}), (\ref{ball}) and (\ref{b-c}), we estimate for any fixed $R>0$
 \bna
 {\nonumber}\int_{\mathbb{B}_R}e^{2\eta_\infty}dx&=&\lim_{k\ra\infty}\int_{\mathbb{B}_{R}(0)}
 \f{\widetilde{u}_k(\widetilde{x}_k+r_k x)\widetilde{f}_k(\widetilde{x}_k+r_k x,\widetilde{u}_k(\widetilde{x}_k+r_k x))}
 {u_k(x_k)f_k(x_k,u_k(x_k))}dx\\
 {\nonumber}&=&\lim_{k\ra\infty}\int_{\mathbb{B}_{Rr_k}(\widetilde{x}_k)}\widetilde{u}_k(x)\widetilde{f}_k(x,\widetilde{u}_k(x))dx
 \\
 {\label{finite}}&\leq&\limsup_{k\ra\infty}\int_{B_{cRr_k}(x_k)}u_kf_k(x,u_k)dv_g\leq C.
 \ena
 It follows that
 $$\int_{\mathbb{R}^2}e^{2\eta_\infty(x)}dx<\infty.$$
 A result of Chen-Li \cite{C-L} implies that
 \be\label{bub}\eta_\infty(x)=-\log(1+|x|^2/4),\quad
 x\in\mathbb{R}^2.\ee
 It follows from (\ref{ball}) that
 $$\int_{\mathbb{B}_{c^{-1}Rr_k}(\widetilde{x}_k)}\widetilde{u}_k\widetilde{f}_k
 (x,\widetilde{u}_k)e^{\psi(x)}dx\leq\int_{B_{Rr_k}(x_k)}u_kf_k(x,u_k)dv_g\leq
 \int_{\mathbb{B}_{cRr_k}(\widetilde{x}_k)}\widetilde{u}_k\widetilde{f}_k
 (x,\widetilde{u}_k)e^{\psi(x)}dx.$$
 In view of (\ref{C-V}) and (\ref{b-c}), we have
 \bna
 \lim_{R\ra\infty}\lim_{k\ra\infty}\int_{\mathbb{B}_{cRr_k}(\widetilde{x}_k)}\widetilde{u}_k\widetilde{f}_k
 (x,\widetilde{u}_k)e^{\psi(x)}dx&=&\lim_{R\ra\infty}\lim_{k\ra\infty}\int_{\mathbb{B}_{c^{-1}Rr_k}(\widetilde{x}_k)}\widetilde{u}_k\widetilde{f}_k
 (x,\widetilde{u}_k)e^{\psi(x)}dx\\
 &=&\lim_{R\ra\infty}\int_{\mathbb{B}_{c^{-1}R}(0)}e^{2\eta_\infty}dx=\int_{\mathbb{R}^2}e^{2\eta_\infty}dx.
 \ena
 Therefore we obtain by (\ref{bub})
  \be\label{eng}
   \lim_{R\ra\infty}\lim_{k\ra\infty}\int_{B_{Rr_k}(x_k)}u_kf_k(x,u_k)dv_g
   =\int_{\mathbb{R}^2}e^{2\eta_\infty}(x)dx=4\pi.
  \ee

 {\it Step 2.} {\it Multi-bubble analysis.}\\

 In this step, we shall prove that there exists some positive integer
 $\ell$ such that the properties $(\mathcal{B}_\ell)$ and
 $(\mathcal{G}_\ell)$ hold. Namely, there exist $\ell$ sequences of points $(x_{i,k})\subset \Sigma$ such that
 $x_{i,k}\ra
 x_i^\ast$ as $k\ra\infty$, $1\leq i\leq \ell$, and the following are satisfied:\\

 \noindent
$(\mathcal{B}_\ell^1)$ For every $i:1\leq i\leq \ell$, letting
$r_{i,k}>0$ be given by (\ref{scal}), $(U_i,\phi_i;\{x^1,x^2\})$ be
an isothermal coordinate system near $x_i^\ast$, where $U_i\subset\Sigma$ is a neighborhood of $x_i^\ast$,
$\phi_i:U_i\ra \Omega_i\subset\mathbb{R}^2$ is a diffeomorphism
with $\phi_i(x_i^\ast)=(0,0)$, and letting $\eta_{i,k}$ be given by
(\ref{blow-up-funct}), we have that $r_{i,k}\ra 0$ as $k\ra\infty$
and
$$\eta_{i,k}(x)\ra \eta_\infty(x)=-\log(1+|x|^2/4)\quad{\rm in}\quad
C^1_{\rm loc}(\mathbb{R}^2)\quad{\rm as}\quad k\ra\infty;$$
$(\mathcal{B}_\ell^2)$ For all $1\leq i\not=j\leq \ell$,
$$\f{d_g(x_{i,k},x_{j,k})}{r_{i,k}}\ra\infty\quad{\rm as}\quad k\ra\infty;$$
 $(\mathcal{B}_\ell^3)$ The following energy identity holds
 $$\lim_{R\ra\infty}\lim_{k\ra\infty}\int_{\cup_{i=1}^\ell{B}_{Rr_{i,k}}(x_{i,k})}
 u_kf_k(x,u_k)dv_g= 4\pi \ell;$$
$(\mathcal{G}_\ell)$ There exists a constant $C>0$  such that
$$R_{\ell,k}^2(x)u_k(x)f_k(x,u_k(x))\leq C$$
for all $x\in \Sigma$ and all $k\in\mathbb{N}$. Here
\be\label{R}R_{\ell,k}(x)=\min_{1\leq i\leq \ell}d_g(x,x_{i,k}).\ee

From Step 1, we know that $(\mathcal{B}_1)$ holds. Suppose for some
$\ell\geq 1$, $(\mathcal{B}_\ell)$ holds but $(\mathcal{G}_\ell)$
does not hold. Choose $x_{\ell+1,k}\in \Sigma$ satisfying
 \bea{\nonumber}
 R^2_{\ell,k}(x_{\ell+1,k})u_k(x_{\ell+1,k})f_k(x_{\ell+1,k},u_k(x_{\ell+1,k}))&=&\max_{x\in\Sigma}
 R_{\ell,k}^2(x)u_k(x)f_k(x,u_k(x))\\&\ra&+\infty\quad{\rm as}\quad
 k\ra\infty.\label{R-0}
 \eea  Let $r_{\ell+1,k}>0$ be as defined in (\ref{scal}). It follows from
 (\ref{scal}), (\ref{R}), and (\ref{R-0}) that $r_{\ell+1,k}\ra 0$ as
 $k\ra\infty$ and
 \be\label{rl}\lim_{k\ra\infty}\f{d_g(x_{\ell+1,k},x_{i,k})}{r_{\ell+1,k}}=+\infty,
 \quad\forall 1\leq i\leq \ell.\ee
 Also we claim that
 \be\label{ik}\lim_{k\ra\infty}\f{d_g(x_{\ell+1,k},x_{i,k})}{r_{i,k}}=+\infty,
 \quad\forall 1\leq i\leq \ell.\ee
 Suppose not. There exists some constant $C$ such that for some $1\leq i\leq
 \ell$, there holds
 $${d_g(x_{\ell+1,k},x_{i,k})}\leq C{r_{i,k}}\quad
 {\rm for\,\,\, all}\quad k.$$ Hence we have
 \be\label{rrr}
 R_{\ell,k}^2(x_{\ell+1,k})u_k(x_{\ell+1,k})f_k(x_{\ell+1,k},u_k(x_{\ell+1,k}))\leq
 Cr_{i,k}^2
 u_k(x_{\ell+1,k})f_k(x_{\ell+1,k},u_k(x_{\ell+1,k}))
 \ee
 By $(\mathcal{B}_\ell^1)$, we estimate
 \bna
 r_{i,k}^2u_k(x_{\ell+1,k})f_k(x_{\ell+1,k},u_k(x_{\ell+1,k}))&=&\f{1+o(1)}{\pi}
 \int_{\mathbb{B}_{r_{i,k}}(\widetilde{x}_{i,k})}\widetilde{u}_k(x)\widetilde{f}_k(x,\widetilde{u}_k(x))
 e^{\psi_i(x)}dx\\
 &\leq&\f{1+o(1)}{\pi}\int_\Sigma u_k(x)f_k(x,u_k(x))dv_g.
 \ena
 This together with (\ref{bdd}) implies that $r_{i,k}^2u_k(x_{\ell+1,k})f_k(x_{\ell+1,k},u_k(x_{\ell+1,k}))$
 is a bounded sequence, and whence  (\ref{rrr}) implies that
 $R_{\ell,k}^2(x_{\ell+1,k})u_k(x_{\ell+1,k})f_k(x_{\ell+1,k},u_k(x_{\ell+1,k}))$
 is bounded. This contradicts (\ref{R-0}). Hence our claim
 (\ref{ik}) holds, and thus $(\mathcal{B}_{\ell+1}^2)$ holds.

 Assume $x_{\ell+1,k}\ra x_{\ell+1}^\ast$ as $k\ra\infty$. Take an isothermal
 coordinate system $(U_{\ell+1},\phi_{\ell+1};\{x^1,x^2\})$ near $x_{\ell+1}^\ast$,
 where $U_{\ell+1}$ is a neighborhood of $x_{\ell+1}^\ast$,
 $\phi_{\ell+1}:U_{\ell+1}\ra\Omega_{\ell+1}\subset\mathbb{R}^2$ is a diffeomorphism
 with $\phi_{\ell+1}(x_{\ell+1}^\ast)=(0,0)$. In
 this coordinate system, the metric $g$ can be represented by
 $$g=e^{\psi_{\ell+1}}(d{x^1}^2+d{x^2}^2)$$
 for some smooth function $\psi_{\ell+1}:\Omega_{\ell+1}\ra\mathbb{R}$ with
 $\psi_{\ell+1}(0,0)=0$. Also we have
 $\nabla_g=e^{-\psi_{\ell+1}}\nabla_{\mathbb{R}^2}$ and
 $\Delta_g=-e^{-\psi_{\ell+1}}\Delta_{\mathbb{R}^2}$.

 Define
 $$v_{\ell+1,k}(x)=\f{\widetilde{u}_k(\widetilde{x}_{\ell+1,k}+r_{\ell+1,k}x)}{u_k(x_{\ell+1,k})}$$
 for $x\in\Omega_{\ell+1,k}=\{x\in\mathbb{R}^2:
 \widetilde{x}_{\ell+1,k}+r_{\ell+1,k}x\in\Omega_{\ell+1}\}$,
 where $\widetilde{x}_{\ell+1,k}=\phi_{\ell+1}(x_{\ell+1,k})$, $\widetilde{u}_k=u_k\circ
 \phi_{\ell+1}^{-1}$.
 Now we prove that
 \be\label{vl1}v_{\ell+1,k}\ra 1\quad{\rm in}\quad C^1_{\rm loc}(\mathbb{R}^2)\quad {\rm as}\quad
 k\ra\infty.\ee
 In view of (\ref{e1}), $v_{\ell+1,k}$ satisfies the equation
 \bea\label{11k}{\nonumber}
 -\Delta_{\mathbb{R}^2}v_{\ell+1,k}(x)&=&e^{\psi_{\ell+1}(\widetilde{x}_{\ell+1,k}+r_{\ell+1,k}x)}\f{\widetilde{f}_k(
 \widetilde{x}_{\ell+1,k}+r_{\ell+1,k}x,
 \widetilde{u}_k(\widetilde{x}_{\ell+1,k}+r_{\ell+1,k}x))}{u_k^2(x_{\ell+1,k})f_k(x_{\ell+1,k},u_k(x_{\ell+1,k}))}\\[1.5ex]
 &&-e^{\psi_{\ell+1}(\widetilde{x}_{\ell+1,k}+r_{\ell+1,k}x)}r_{\ell+1,k}^2\widetilde{\tau}_k
 (\widetilde{x}_{\ell+1,k}+r_{\ell+1,k}x) v_{\ell+1,k}(x)\eea on $\Omega_{\ell+1,k}$,
 where $\widetilde{f}_k(x,t)=f_k(\phi_{\ell+1}^{-1}(x),t)$.
 By (\ref{R-0}), we have
 \bea\nonumber
 &&\widetilde{R}_{\ell,k}^2(\widetilde{x}_{\ell+1,k}+r_{\ell+1,k}x)\widetilde{u}_k(\widetilde{x}_{\ell+1,k}+r_{\ell+1,k}x)
 \widetilde{f}_k(\widetilde{x}_{\ell+1,k}+r_{\ell+1,k}x,\widetilde{u}_k(\widetilde{x}_{\ell+1,k}+r_{\ell+1,k}x))\\
 \label{Rleq}&&\leq
 R_{\ell,k}^2(x_{\ell+1,k})u_k(x_{\ell+1,k})f_k(x_{\ell+1,k},u_k(x_{\ell+1,k})),
 \eea
 where $\widetilde{R}_{\ell,k}={R}_{\ell,k}\circ \phi_{\ell+1}^{-1}$.
 Fix any $i$, $1\leq i\leq \ell$. If $x_{\ell+1}^\ast\not=x_i^\ast$,
 noting that $d_g(\phi_{\ell+1}^{-1}(\widetilde{x}_{\ell+1,k}+r_{\ell+1,k}x),x_{i,k})\ra d_g(x_{\ell+1}^\ast,x_i^\ast)$
 and $d_g(x_{\ell+1,k},x_{i,k})\ra d_g(x_{\ell+1}^\ast,x_i^\ast)$ as
 $k\ra\infty$, we then have \be\label{ooo}d_g(\phi_{\ell+1}^{-1}(\widetilde{x}_{\ell+1,k}+r_{\ell+1,k}x),x_{i,k})=(1+o(1))
 d_g(x_{\ell+1,k},x_{i,k}),\ee where $o(1)\ra 0$ as $k\ra\infty$ uniformly in $x\in\mathbb{B}_R(0)$.
 If $x_{\ell+1}^\ast=x_i^\ast$, since the Riemannian distance and the Euclidean distance
 are equivalent in the same local coordinate system, we then have $|\phi_{\ell+1}(x_{\ell+1,k})-\phi_{\ell+1}(x_{i,k})|
 =(1+o(1))d_g(x_{\ell+1,k},x_{i,k})$. Recalling (\ref{rl}), we
 obtain for all $x\in\mathbb{B}_R(0)$
 \bna
 d_g(\phi_{\ell+1}^{-1}(\widetilde{x}_{\ell+1,k}+r_{\ell+1,k}x),x_{i,k})&=&
 (1+o(1))|\widetilde{x}_{\ell+1,k}+r_{\ell+1,k}x-\phi_{\ell+1}(x_{i,k})|\\
 &=&(1+o(1))d_g(x_{\ell+1,k},x_{i,k}).
 \ena
 Hence we have (\ref{ooo}) in any case. Combining (\ref{Rleq}) and (\ref{ooo}), we obtain for  $x\in
 \mathbb{B}_R(0)$
 \bea\label{ineq}
 &&{v}_{\ell+1,k}(x)\f{\widetilde{f}_k(\widetilde{x}_{\ell+1,k}+r_{\ell+1,k}x,\widetilde{u}_k(\widetilde{x}_{\ell+1,k}+r_{\ell+1,k}x))}
 {f_k(x_{\ell+1,k},u_k(x_{\ell+1,k}))}\nonumber\\
 &&\leq\f{\inf_{1\leq i\leq \ell}d_g(x_{\ell+1,k},x_{i,k})^2}{\inf_{1\leq i\leq
 \ell}d_g(\phi_{\ell+1}^{-1}(\widetilde{x}_{\ell+1,k}+r_{\ell+1,k}x),x_{i,k})^2}=1+o(1),
 \eea
 where $o(1)\ra 0$ uniformly in $x\in\mathbb{B}_R(0)$. From (H4), we know that there exists $t_0>0$ such that
 \be\label{f-geq}\f{f_k(x,t_2)}{f_k(x,t_1)}\geq e^{t_2^2-t_1^2}\,\, {\rm for\,\,all}\,\, t_1, t_2\geq t_0,
 \,\,{\rm and\,\,all}\,\, x\in\Sigma.\ee
 If there exist some $R_0>0$ and a sequence of points $(z_k)\subset\mathbb{B}_{R_0}(0)$ such that
  $v_{\ell+1,k}(z_k)\ra \alpha>1$ as $k\ra\infty$, then we conclude by
 (\ref{f-geq}) and (H5) that
 $${v}_{\ell+1,k}(z_k)\f{\widetilde{f}_k(\widetilde{x}_{\ell+1,k}+r_{\ell+1,k}z_k,
 \widetilde{u}_k(\widetilde{x}_{\ell+1,k}+r_{\ell+1,k}z_k))}{f_k(x_{\ell+1,k},u_k(x_{\ell+1,k}))}
 \geq \f{\alpha+1}{2}>1$$ for sufficiently large $k$,
 which contradicts (\ref{ineq}). Therefore we obtain
 $$\limsup_{k\ra\infty}\|v_{\ell+1,k}\|_{L^\infty(\mathbb{B}_R(0))}\leq 1,\quad\forall R>0.$$
 When $v_{\ell+1,k}(x)>1$, we
 have by (\ref{11k}) and (\ref{ineq}),
 $\Delta_{\mathbb{R}^2}v_{\ell+1,k}(x)=o(1)$,
 where $o(1)$ is the same meaning as that of (\ref{ineq}).
 When $v_{\ell+1,k}(x)\leq 1$,
 using (H4) and (H5), we also have
 $\Delta v_{\ell+1,k}(x)=o(1)$,
 where $o(1)\ra 0$ as $k\ra\infty$ uniformly in all $x$ satisfying $v_{\ell+1,k}(x)\leq
 1$ for sufficiently large $k$. Now applying elliptic estimates to
 equation (\ref{11k}), we obtain
 $$v_{\ell+1,k}\ra v_{\ell+1,\infty}\quad{\rm in}\quad C^1_{\rm loc}(\mathbb{R}^2),$$
 where $v_{\ell+1,\infty}$ is a solution to
 $$\le\{\begin{array}{lll}
 -\Delta_{\mathbb{R}^2}v_{\ell+1,\infty}=0\quad{\rm in}\quad\mathbb{R}^2\\[1.5ex]
 0\leq v_{\ell+1,\infty}\leq 1.
 \end{array}\ri.$$
 Note that $v_{\ell+1,\infty}(0)=1$. The Liouville theorem for harmonic functions leads to $v_{\ell+1,\infty}\equiv 1$.
 Whence (\ref{vl1}) holds.

 Define another sequence of blow-up functions by
 \be\label{blow-up-l}\eta_{\ell+1,k}(x)=u_k(x_{\ell+1,k})(\widetilde{u}_k(\widetilde{x}_{\ell+1,k}
 +r_{\ell+1,k}x)-u_k(x_{\ell+1,k})),\quad x\in\Omega_{\ell+1,k}.\ee
 In the following, we will prove that $(\mathcal{B}_{\ell+1}^1)$ and
 $(\mathcal{B}_{\ell+1}^3)$ hold. By (\ref{e1}), $\eta_{\ell+1,k}$ satisfies the
 equation
 \bea{\nonumber}-\Delta_{\mathbb{R}^2}\eta_{\ell+1,k}(x)&=&e^{\psi_{\ell+1}(\widetilde{x}_{\ell+1,k}+r_{\ell+1,k}x)}
 \f{\widetilde{f}_k(\widetilde{x}_{\ell+1,k}+r_{\ell+1,k}x,
 \widetilde{u}_k(\widetilde{x}_{\ell+1,k}+r_{\ell+1,k}x))}{f_k(x_{\ell+1,k},u_k(x_{\ell+1,k}))}\\[1.5ex]&&
 -e^{\psi_{\ell+1}(\widetilde{x}_{\ell+1,k}+r_{\ell+1,k}x)}\widetilde{\tau}_k(\widetilde{x}_{\ell+1,k}+r_{\ell+1,k}x)
 r_{\ell+1,k}^2u_k^2(x_{\ell+1,k})v_{\ell+1,k}(x)\label{eta-l1}\eea
 on $\Omega_{\ell+1,k}$.
 We claim that for any fixed $R>0$,
 \be\label{upper}\limsup_{k\ra\infty}\eta_{\ell+1,k}(x)\leq 0
 \,\,{\rm uniformly\,\,in}\,\,x\in\mathbb{B}_R(0).\ee
 For otherwise, we may take a sequence of points
 $(y_k)\subset\mathbb{B}_R(0)$ such that $\eta_{\ell+1,k}(y_k)\geq
 \beta>0$ for all sufficiently large $k$. By (H4), (H5) and
 (\ref{vl1}), we obtain
 \bna
 \f{\widetilde{f}_k(\widetilde{x}_{\ell+1,k}+r_{\ell+1,k}y_k,
 \widetilde{u}_k(\widetilde{x}_{\ell+1,k}+r_{\ell+1,k}y_k))}{f_k(x_{\ell+1,k},u_k(x_{\ell+1,k}))}&=&
 (1+o(1))e^{\widetilde{u}_k^2(\widetilde{x}_{\ell+1,k}+r_{\ell+1,k}y_k)-u_k^2(x_{\ell+1,k})}\\
 &=&(1+o(1))e^{(2+o(1))\eta_{\ell+1,k}(y_k)}\\
 &\geq& 1+2\beta+o(1).
 \ena
 This together with (\ref{ineq}) leads to
 $$1+2\beta+o(1)\leq 1+o(1),$$
 which is impossible when $k$ is sufficiently large. Hence our claim
 (\ref{upper}) holds. By (\ref{vl1}), using the same method of deriving
 (\ref{rk-0}), we conclude
 \be\label{rl10}r_{\ell+1,k}^2u_{k}^2(x_{\ell+1,k})\ra 0\quad{\rm as}\quad k\ra\infty.\ee
  Combining (\ref{vl1}) and (\ref{blow-up-l})-(\ref{rl10}),
  similarly as we did in Step 1, we arrive at
  $$\eta_{\ell+1,k}(x)\ra \eta_\infty(x)\,\, {\rm in}\,\,
  C^1_{\rm loc}(\mathbb{R}^2)\,\,{\rm as}\,\, k\ra\infty,$$
  where $\eta_\infty(x)=-\log(1+|x|^2/4)$ is the unique solution to
  (\ref{bubble}). Hence $(\mathcal{B}_{\ell+1}^1)$ holds.

  Moreover, using the same method for proving (\ref{eng}), we arrive
  at
  \bna
  \lim_{R\ra\infty}\lim_{k\ra\infty}\int_{B_{Rr_{\ell+1,k}}(x_{\ell+1,k})}u_kf_k(x,u_k)dv_g
  =\int_{\mathbb{R}^2}e^{2\eta_\infty(x)}dx=4\pi.
  \ena
  Thus $(\mathcal{B}_{\ell+1}^3)$ holds.

  Actually, we have proved that if $(\mathcal{B}_\ell)$ holds but $(\mathcal{G}_\ell)$ does not hold,
  then $(\mathcal{B}_{\ell+1})$ holds. Note that
  \be\label{stops}\int_{\Sigma}u_kf_k(x,u_k)dv_g\geq \sum_{i=1}^{\ell+1}
  \int_{B_{Rr_{i,k}}(x_{i,k})}u_kf_k(x,u_k)dv_g=4(\ell+1)\pi.\ee
  In view of (\ref{bdd}), the process must be terminate after finite steps. This ends the proof of Step 2.\\

  {\it Step 3.} {\it Exhaustion of blow-up points.}\\

  It follows from Step 2 that there exists some
  $\ell\in\mathbb{N}\setminus\{0\}$ and $\ell$ sequences of points $(x_{i,k})$, $i=1,\cdots,\ell$,
  such that $(\mathcal{B}_\ell)$ and $(\mathcal{G}_\ell)$ hold. If there
  exists a sequence of points $(x_{\ell+1,k})$ of $\Sigma$ such that
  after extracting a new subsequence from the previous one, $(\mathcal{B}_{\ell+1})$ and
  $(\mathcal{G}_{\ell+1})$ hold, we add this sequence of points, and so
  on. The process necessarily terminates because of (\ref{bdd})
  and (\ref{stops}).
  Therefore there exists some $N\in\mathbb{N}\setminus\{0\}$ and $N$
  sequences of points $(x_{i,k})$, $i=1,\cdots,N$, such that
  $(\mathcal{B}_N)$ and $(\mathcal{G}_N)$ hold and such that, given
  any sequence of points $(x_{N+1,k})$, it is impossible to extract
  a new subsequence from the previous one such that $(\mathcal{B}_{N+1})$
  and $(\mathcal{G}_{N+1})$ hold with sequences $(x_{i,k})$,
  $i=1,\cdots,N+1$.\\

  {\it Step 4.} {\it Convergence away from blow-up points.}\\

  Set $\mathcal{S}=\{x_1^\ast,\cdots,x_N^\ast\}$. We will prove that
  $u_k\ra u_\infty$ in $C^1_{\rm loc}(\Sigma\setminus\mathcal{S})$.
  In view of $(\mathcal{G}_N)$, given any compact set $K\subset\Sigma\setminus\mathcal{S}$,
  there exists a constant $C$ such that
  $$u_k(x)f_k(x,u_k(x))\leq C\,\,{\rm for\,\,all}\,\, x\in K
  \,\,{\rm and\,\,all}\,\, k.$$
  If $u_k(x)>1$ for some $x\in K$, then $f_k(x,u_k(x))\leq C_K$. If
  $u_k(x)\leq 1$ for some $x\in K$, then (H2) implies that
  $f_k(x,u_k(x))$ is  bounded uniformly in $x$ with $u_k(x)\leq 1$. Thus, for all $x\in
  K$, $f_k(x,u_k(x))$ is bounded in $L^\infty(K)$. In view of (\ref{tauk}) and (\ref{Lp-bdd}), applying elliptic
  estimates to the equation
  $$\Delta_gu_k(x)+\tau_k(x)u_k(x)=f_k(x,u_k(x)),\quad {x\in K},$$
  we obtain the convergence $u_k\ra u_0$ in $C^1_{\rm
  loc}(\Sigma\setminus\mathcal{S})$.\\

  Combining the above four steps, we complete the proof of
  Proposition 3.1. $\hfill\Box$

 \section{Gradient estimate}

 Let $u_k\geq 0$ be a sequence of solutions to (\ref{e1}). In this section we shall
 establish a gradient estimate on $u_k$, which can be viewed as a version on manifolds of
 (\cite{Druet}, Proposition 2). Precisely we have the following result.\\

 \noindent{\bf Proposition 4.1} {\it Let $(\Sigma,g)$ be a compact Riemannian surface without boundary,
$f_k$ be a sequence of functions satisfying {\rm (H1)-(H5)}, and
$u_k\geq 0$ be a sequence of smooth solutions to equation
(\ref{e1}) such that (\ref{beata}) holds. Assume that
$\max_\Sigma u_k\ra+\infty$ as $k\ra\infty$. Let
$N\in\mathbb{N}\setminus\{0\}$ and the sequences $x_{i,k}$,
$i=1,\cdots,N$, be given by Proposition 3.1. Then there exists a
uniform constant $C$ such that
$$R_{N,k}(x)u_k(x)|\nabla_gu_k(x)|\leq C$$
for all $x\in\Sigma$ and all $k$, where $R_{N,k}(x)$ is defined as in (\ref{R}).}\\

\noindent{\it Proof.} Choose $y_k\in\Sigma$ such that
\be\label{41}R_{N,k}(y_k)u_k(y_k)|\nabla_gu_k(y_k)|=\max_{x\in\Sigma}R_{N,k}(x)
u_k(x)|\nabla_gu_k(x)|.\ee Suppose by contradiction that
\be\label{contrad}R_{N,k}(y_k)u_k(y_k)|\nabla_gu_k(y_k)|\ra+\infty\quad{\rm
as}\quad k\ra\infty.\ee Set \be\label{sk}s_k=R_{N,k}(y_k).\ee By
Proposition 3.1, we have $u_k\ra u_\infty$ in $C^1_{\rm
loc}(\Sigma\setminus\{x_1^\ast,\cdots,x_N^\ast\})$, which together
with (\ref{contrad}) implies that $s_k\ra 0$ as $k\ra\infty$.
Without loss of generality, we may assume that $y_k\ra x_1^\ast$ as
$k\ra\infty$, $x_1^\ast=\cdots=x_\ell^\ast$ for some $1\leq \ell\leq
N$, and $x_j^\ast\not=x_1^\ast$ for any $j\in\{\ell+1,\cdots,N\}$.
Take an isothermal coordinate system $(U,\phi;\{x^1,x^2\})$ near
$x_1^\ast$, where $U$ is a neighborhood of $x_1^\ast\in\Sigma$,
$\phi: U\ra\Omega\subset\mathbb{R}^2$ is a diffeomorphism with
$\phi(x_1^\ast)=(0,0)$. In this coordinate system the
metric $g$ can be represented by $g=e^{\psi}(d{x^1}^2+d{x^2}^2)$,
where $\psi:\Omega\ra\mathbb{R}$ is a smooth function with
$\psi(0,0)=0$. Denote $\widetilde{y}_k=\phi(y_k)$,
$\widetilde{u}_k=u_k\circ\phi^{-1}$. We
set
$$v_k(y)=\widetilde{u}_k(\widetilde{y}_k+s_ky)$$
for $y\in\Omega_k=\{y\in\mathbb{R}^2:
\widetilde{y}_k+s_ky\in\Omega\}$. Define
$$y_{i,k}=\f{\widetilde{x}_{i,k}-\widetilde{y}_k}{s_k}\in \Omega_k,\quad i=1,\cdots,\ell,$$
and
$$\widetilde{S}_k=\{y_{1,k},\cdots, y_{\ell,k}\}.$$
Since $s_k\ra 0$, we have $\Omega_k\ra\mathbb{R}^2$ as $k\ra\infty$.
Denote
$$\widetilde{S}=\lim_{k\ra\infty}\widetilde{S}_k.$$
By (\ref{sk}) and the fact $\psi(0,0)=0$, we have
 \bna
 d_{\mathbb{R}^2}(0,\widetilde{S}_k)&=&\inf_{1\leq i\leq
 \ell}|y_{i,k}|=\inf_{1\leq i\leq
 \ell}\f{|\widetilde{x}_{i,k}-\widetilde{y}_k|}{s_k}\\
 &=&\inf_{1\leq i\leq
 \ell}\f{(1+o(1))d_g(x_{i,k},y_k)}{s_k}\\
 &=&1+o(1),
 \ena
and thus
\be\label{4.3}{d}_{\mathbb{R}^2}(0,\widetilde{S})=1,\ee
where ${ d}_{\mathbb{R}^2}(\cdot,\cdot)$ denotes the Euclidean
distance of $\mathbb{R}^2$. Clearly, $v_k(y)$ satisfies
 \be\label{vvv}
 -\Delta_{\mathbb{R}^2}v_k(y)=e^{\psi(\widetilde{y}_k+s_ky)}s_k^2\le(\widetilde{f}_k(\widetilde{y}_k+s_ky,
 \widetilde{u}_k(\widetilde{y}_k+s_ky))-
 \widetilde{\tau}_k(\widetilde{y}_k+s_ky) v_k(y)\ri)\ee
 for $y\in \Omega_k$.
 By $(iii)$ of Proposition 3.1, we have
 \be\label{point}\widetilde{R}_{N,k}(\widetilde{y}_k+s_ky)^2v_k(y)
 \widetilde{f}_k(\widetilde{y}_k+s_ky,v_k(y))\leq C\ee
 for some constant $C$ independent of $k$. Note that
 \bea\nonumber\widetilde{R}_{N,k}(\widetilde{y}_k+s_ky)&=&R_{N,k}(\phi^{-1}(\widetilde{y}_k+s_ky))\\
 \nonumber&=&\inf_{1\leq
 i\leq\ell}d_g(\phi^{-1}(\widetilde{y}_k+s_ky),x_{i,k})\\
 \nonumber&=&(1+o(1))\inf_{1\leq i\leq\ell}d_{\mathbb{R}^2}(\widetilde{y}_k+s_ky,\widetilde{x}_{i,k})\\
 &=&(1+o(1))s_k d_{\mathbb{R}^2}(y,\widetilde{S}_k).\label{RNK}
 \eea
 Combining (\ref{point}) and (\ref{RNK}), we have
 \be\label{vdelta}s_k^2v_k(y)\widetilde{f}_k(\widetilde{y}_k+s_ky,v_k(y))
 \leq \f{C}{d_{\mathbb{R}^2}(y,\widetilde{S}_k)^2},\ee
 which together with (H1) and (H2) leads to
 \be\label{fk}0\leq s_k^2\widetilde{f}_k(\widetilde{y}_k+s_ky,v_k(y))\leq
 \f{C}{d_{\mathbb{R}^2}(y,\widetilde{S}_k)^2}.\ee
  In view of (\ref{Lp-bdd}), we estimate for any $p>1$ and any $R>0$,
 \bea
 {\nonumber}\int_{\mathbb{B}_R(0)}(s_k^2v_k(y))^pdy&=&
 s_k^{2p}\int_{\mathbb{B}_R(0)}\widetilde{u}_k(\widetilde{y}_k+s_ky)^pdy\\
 {\nonumber}&\leq&
 Cs_k^{2p-2}\int_{\Sigma}u_k^pdv_g\\{\label{skvk}}
 &\ra& 0\quad {\rm as}\quad k\ra\infty.
 \eea
 Denote for any $R>0$
 $$A_{R}=\mathbb{B}_R(0)\setminus\cup_{y\in \mathcal{S}}{\mathbb{B}_{1/R}(y)}.$$
 Clearly there exists some $R_0>0$ such that $A_{R/4}$ is necessarily smooth
 bounded domain provided that $R\geq R_0$. Now we take $R\geq R_0$.
 In view of (\ref{tauk}), (\ref{vvv}), (\ref{fk}), and (\ref{skvk}), we arrive at
 $$\lim_{k\ra\infty}\|\Delta_{\mathbb{R}^2}v_k\|_{L^p(A_R)}= 0,\quad\forall R\geq R_0,\quad \forall p>1.$$
 Let $w_k$ satisfy
 $$\le\{\begin{array}{lll}
 -\Delta_{\mathbb{R}^2}w_k=-\Delta_{\mathbb{R}^2}v_k\quad{\rm in}\quad A_R\\[1.5ex]
 w_k=0\quad {\rm on}\quad \p A_R.
 \end{array}\ri.$$
 It follows from (\ref{skvk}) and elliptic estimates that there exists some function $w$
 such that
 $$w_k\ra w\quad{\rm in}\quad C^1(\overline{A_R}).$$
 In particular, $w_k$ is uniformly bounded in $A_R$.
 While $v_k-w_k$ satisfies
 \be\label{v-wk}\le\{\begin{array}{lll}
 -\Delta_{\mathbb{R}^2}(v_k-w_k)=0\quad{\rm in}\quad A_R\\[1.5ex]
 v_k-w_k=v_k\quad {\rm on}\quad \p A_R.
 \end{array}\ri.\ee
 We claim that
 \be\label{vk-0}v_k(0)\ra+\infty\quad{\rm as}\quad k\ra\infty.\ee
 For otherwise, $(v_k(0)-w_k(0))$ would be a bounded sequence. Noting that
 $v_k-w_k$ has a lower bound in $A_{R}$, applying Harnack's
 inequality to (\ref{v-wk}), we obtain
 $$\|v_k-w_k\|_{L^\infty(A_{R/2})}\leq C$$
 for some constant $C$ depending only on $R$, and whence $v_k$ is bounded
 in $C^1(A_{R/4})$. In view of (\ref{4.3}), this leads to
  $$v_k(0)|\nabla_{\mathbb{R}^2}v_k(0)|\leq C.$$
  While (\ref{41}) and (\ref{contrad}) implies
  \be\label{vk0-infty}v_k(0)|\nabla_{\mathbb{R}^2}v_k(0)|\ra+\infty\quad{\rm as}\quad k\ra\infty.\ee
  This is a contradiction. Hence our claim (\ref{vk-0}) follows.

  Replacing $v_k$ by $v_k/v_k(0)$ in the above estimates, we obtain
  \be\label{con-1}\f{v_k}{v_k(0)}\ra 1\quad{\rm in}\quad
  C^1_{\rm loc}(\mathbb{R}^2\setminus\mathcal{S})\ee
  as $k\ra\infty$. For $y\in \Omega_k$, we set
  $$\widetilde{v}_k(y)=\f{v_k(y)-v_k(0)}{|\nabla_{\mathbb{R}^2}v_k(0)|}.$$
  It follows from (\ref{41}) and (\ref{RNK}) that
  $$v_k(y)|\nabla_{\mathbb{R}^2}v_k(y)|\leq (1+o(1))\f{v_k(0)|\nabla_{\mathbb{R}^2}v_k(0)|}
  {d_{\mathbb{R}^2}(y,\widetilde{S}_k)},\quad y\in\Omega_k\setminus \widetilde{S}_k.$$
  This together with (\ref{con-1}) gives
  \be\label{grad}|\nabla_{\mathbb{R}^2}\widetilde{v}_k(y)|\leq \f{1+o(1)}
  {d_{\mathbb{R}^2}(y,\mathcal{S})},\ee
  where $o(1)\ra 0$ as $k\ra\infty$ locally uniformly in
  $y\in\mathbb{R}^2\setminus\mathcal{S}$. Since
  $\widetilde{v}_k(0)=0$, it follows from (\ref{grad}) that
  $\widetilde{v}_k$ is uniformly bounded in $C^1(A_R)$ for any
  $R>0$. In view of (\ref{vvv}) and (\ref{con-1}), we have
  \bea
  {\nonumber}-\Delta_{\mathbb{R}^2}\widetilde{v}_k(y)&=&-(1+o(1))\f{v_k(y)\Delta_{\mathbb{R}^2}v_k(y)}{v_k(0)|\nabla_{\mathbb{R}^2}v_k(0)|}\\
  &=&\f{1+o(1)}{v_k(0)|\nabla_{\mathbb{R}^2}v_k(0)|}e^{\psi(\widetilde{y}_k+s_ky)}s_k^2v_k(y)\le\{
  \widetilde{f}_k(\widetilde{y}_k+s_ky,
 v_k(y))-\widetilde{\tau}_k(\widetilde{y}_k+s_ky) v_k(y)\ri\}\quad\quad{\label{delt-v}}
  \eea
  for $y\in \Omega_k$.
  Similarly to (\ref{skvk}), $s_k^2v_k^2$ is bounded in $L^p_{\rm loc}(\mathbb{R}^2)$ for
  any $p>1$. In view of (\ref{vdelta}) and (\ref{vk0-infty}), applying elliptic estimates to
  the equation (\ref{delt-v}), we have
  \be\label{vk-con}\widetilde{v}_k\ra \widetilde{v}\quad{\rm in}\quad C^1_{\rm loc}(\mathbb{R}^2\setminus \mathcal{S})
  \quad{\rm as}\quad k\ra\infty,\ee
  where $\widetilde{v}$ satisfies
  \be\label{420} \Delta_{\mathbb{R}^2}\widetilde{v}=0\quad{\rm in}\quad \mathbb{R}^2\setminus\mathcal{S},\quad
  \widetilde{v}(0)=0,\quad |\nabla_{\mathbb{R}^2}\widetilde{v}(0)|=1,\ee
  and
  \be\label{gr-vt}|\nabla_{\mathbb{R}^2}\widetilde{v}(y)|\leq \f{1}{d_{\mathbb{R}^2}(y,\mathcal{S})},
  \quad y\in\mathbb{R}^2\setminus\mathcal{S}.\ee
  Let $\hat{y}\in\mathcal{S}$. For any
  $0<r<d_{\mathbb{R}^2}(\hat{y},\mathcal{S}\setminus\{\hat{y}\})/2$, since
  \bna
  \int_{\mathbb{B}_r(\hat{y})}v_k\Delta_{\mathbb{R}^2}v_kdy&=&\int_{\mathbb{B}_r(\hat{y})}
  \widetilde{u}_k(\widetilde{y}_k+s_ky)s_k^2\Delta_{\mathbb{R}^2}\widetilde{u}_k(\widetilde{y}_k+s_ky)dy\\
  &=&\int_{\mathbb{B}_{s_kr}(\widetilde{y}_k+s_k\hat{y})}\widetilde{u}_k(x)\Delta_{\mathbb{R}^2}
  \widetilde{u}_k(x)dx\\
  &=&-\int_{\phi^{-1}(\mathbb{B}_{s_kr}(\widetilde{y}_k+s_k\hat{y}))}u_k\Delta_g
  u_kdv_g,
  \ena
  we get by (\ref{tauk}), (\ref{bdd}) and (\ref{Lp-bdd})
  $$\le|\int_{\mathbb{B}_{r}(\hat{y})}v_k\Delta_{\mathbb{R}^2}v_kdy\ri|\leq \int_\Sigma
  \le(u_kf_k(x,u_k)+\tau_k u_k^2\ri)dv_g\leq C.$$
  Similarly we have by (\ref{energy-bounded}) $$\int_{\mathbb{B}_r(\hat{y})}
  |\nabla_{\mathbb{R}^2}v_k|^2dy\leq
  \int_\Sigma|\nabla_gu_k|^2dv_g\leq C.$$
  It then follows that
  $$\int_{\p\mathbb{B}_r(\hat{y})}v_k\p_\nu v_k d\sigma=\int_{\mathbb{B}_r(\hat{y})}
  |\nabla_{\mathbb{R}^2}v_k|^2dy-\int_{\mathbb{B}_r(\hat{y})}v_k\Delta_{\mathbb{R}^2}v_kdy=O(1).$$
  While (\ref{con-1}) and (\ref{vk-con}) lead to
  $$\int_{\p\mathbb{B}_r(\hat{y})}v_k\p_\nu v_k d\sigma=v_k(0)|\nabla_{\mathbb{R}^2}v_k(0)|
  \le(\int_{\p\mathbb{B}_r(\hat{y})}\p_{\nu}\widetilde{v}d\sigma+o(1)\ri).$$
  This together with (\ref{vk0-infty}) gives for any $0<r<d_{\mathbb{R}^2}(\hat{y},\mathcal{S}\setminus\{\hat{y}\})/2$
  $$\int_{\p\mathbb{B}_r(\hat{y})}\p_\nu \widetilde{v}d\sigma=0,$$
  which leads to
  $$
  \f{d}{dr}\le(\f{1}{2\pi
  r}\int_{\p\mathbb{B}_r(\hat{y})}\widetilde{v}d\sigma\ri)=
  \f{1}{2\pi
  r}\int_{\p\mathbb{B}_r(\hat{y})}\p_\nu\widetilde{v}d\sigma=0.
  $$
  Hence there exists some constant $\alpha$ depending only on
  $\hat{y}$ such that
  \be\label{mea-const}\f{1}{2\pi
  r}\int_{\p\mathbb{B}_{r}(\hat{y})}\widetilde{v}d\sigma=\alpha,\quad\forall
  0<r<d_{\mathbb{R}^2}(\hat{y},\mathcal{S}\setminus\{\hat{y}\})/2.\ee
  Given any $y\in\p \mathbb{B}_{r}(\hat{y})$. (\ref{mea-const})
  permits us to take $y^\ast\in \p \mathbb{B}_{r}(\hat{y})$ such
  that $\widetilde{v}(y^\ast)=\alpha$. It then follows from
  (\ref{gr-vt}) that $|\widetilde{v}(y)-\alpha|\leq \pi$. This
  indicates that $\widetilde{v}$ is bounded near $\hat{y}$. Since
  this is true for all $\hat{y}\in\mathcal{S}$, we conclude that
  $\widetilde{v}$ is a smooth harmonic function in $\mathbb{R}^2$.
  By the mean value equality,
  $$\int_{\p\mathbb{B}_R(0)}\widetilde{v}d\sigma=0,\quad \forall R>0.$$
  This together with (\ref{gr-vt}) implies that $\widetilde{v}$ is
  bounded in $L^\infty(\mathbb{R}^2)$. Actually we can take $z\in\p\mathbb{B}_R(0)$
  such that $\widetilde{v}_k(z)=0$, in view of (\ref{gr-vt}), we then have
  for all $y\in\p\mathbb{B}_R(0)$
  \bna
  |\widetilde{v}_k(y)|=|\widetilde{v}_k(y)-\widetilde{v}_k(z)|\leq
  \pi R\sup_{\p\mathbb{B}_R(0)}|\nabla_{\mathbb{R}^2}
  \widetilde{v}|\leq 2\pi,
  \ena
  provided that $R>2\sup_{\hat{y}\in\mathcal{S}}|\hat{y}|$.
   Note again that
  $\widetilde{v}(0)=0$. Applying the Liouville theorem to (\ref{420}), we have $\widetilde{v}\equiv
  0$, which contradicts the fact that
  $|\nabla_{\mathbb{R}^2}\widetilde{v}(0)|=1$. This completes the proof of the
  proposition. $\hfill\Box$

  \section{Quantization}

  In this section we prove quantization results for equation
  (\ref{e1}). Let $x_1^\ast,\cdots,x_N^\ast$ be as in Proposition 3.1.
  For some $1\leq i\leq N$, $x_i^\ast$ is called a {\it simple} blow-up point
  if $N=1$ or $x_j\not=x_i$ for all $j\in\{1,\cdots,N\}\setminus\{i\}$; Otherwise
  we call $x_i^\ast$ a {\it non-simple} blow-up point. In the following, we distinguish between
  these two types of points to proceed.

  \subsection{Quantization for simple blow-up points\\}

  Let $x_i^\ast$ be a simple blow-up point. Take an isothermal coordinate system $(U_i,\phi_i;\{x^1,x^2\})$ near $x_i^\ast$,
  where $U_i\subset\Sigma$ is a neighborhood  of $x_i^\ast$ such that $x_j^\ast\not\in \overline{U}_i$,
  the closure of $U_i$, for all $j\in\{1,\cdots,N\}\setminus\{i\}$. As before $\phi_i:U_i\ra\Omega\subset\mathbb{R}^2$ is a diffeomorphism
  with $\phi_i(x_i^\ast)=(0,0)$. Particularly we can find some $\delta>0$ such that $\mathbb{B}_{2\delta}(0)\subset
  \Omega$. In this coordinate system, the metric $g$ writes as $g=e^{\psi_i}(d{x^1}^2+d{x_2}^2)$
  for some smooth function $\psi_i:\Omega\ra\mathbb{R}$
  with $\psi_i(0,0)=0$. In this subsection we prove the following quantization result. \\

  \noindent{\bf Proposition 5.1} {\it Let $u_k$, $u_\infty$, $\tau_k$, $\tau_\infty$, $x_{i,k}$ and $x_i^\ast$ be as in Proposition 3.1.
  Suppose that $x_i^\ast$
  is a simple blow-up point. Then up to a subsequence, there exists some positive integer $I^{(i)}$
  such that
  \be\label{identity}
  \lim_{k\ra\infty}\int_{U_i}
  (|\nabla_gu_k|^2+\tau_ku_k^2)dv_g=
  \int_{U_i}(|\nabla_gu_\infty|^2+\tau_\infty
  u_\infty^2)dv_g+4\pi I^{(i)},
  \ee
  where $U_i$ is a neighborhood of $x_i^\ast$ as above.}\\

  In the coordinate system $(U_i,\phi_i;\{x^1,x^2\})$, we write
  $\widetilde{x}_{i,k}=\phi_i^{-1}(x_{i,k})$,
  $\widetilde{u}_k(x)=u_k(\phi_i^{-1}(x))$,
  $\widetilde{\tau}_k(x)=\tau_k(\phi_i^{-1}(x))$ and $\widetilde{f}_k(x,\widetilde{u}_k(x))=f_k(\phi_i^{-1}(x),u_k(\phi_i^{-1}(x)))$
  for any $x\in\Omega$. Moreover for $0<s<t<\delta$ we define the spherical mean of $\widetilde{u}_k$, the total energy and
  the neck energy of $\widetilde{u}_k$ around $\widetilde{x}_{i,k}$ by
  \be\label{mean-val}\varphi_k(t)=\varphi_k^{(i)}(t)=\f{1}{2\pi t}\int_{\p\mathbb{B}_t(\widetilde{x}_{i,k})}\widetilde{u}_kd\sigma,\ee
  \be\label{energy}\Lambda_{k}(t)=\Lambda_k^{(i)}(t)=\int_{\mathbb{B}_t(\widetilde{x}_{i,k})}
  \widetilde{u}_k\widetilde{f}_k(x,\widetilde{u}_k)dx,\ee
  and
  \be\label{neck}N_{k}(s,t)=N_k^{(i)}(s,t)=\int_{\mathbb{B}_t(\widetilde{x}_{i,k})\setminus \mathbb{B}_s(\widetilde{x}_{i,k})}
  \widetilde{u}_k\widetilde{f}_k(x,\widetilde{u}_k)dx\ee
  respectively.
  We say that the property $(\mathcal{H}_\ell)$ holds if there
  exist sequences
  $$s_k^{(0)}=0<r_k^{(1)}<s_k^{(1)}<\cdots<r_k^{(\ell)}<s_k^{(\ell)}=o(1)$$
  such that the following hypotheses are satisfied:\\

  \noindent$(\mathcal{H}_{\ell,1})$\,\,\, $\lim\limits_{k\ra\infty}{r_k^{(j)}}/{s_k^{(j)}}=\lim\limits_{k\ra\infty}
  {s_k^{(j-1)}}/{r_k^{(j)}}=0$ for all $1\leq j\leq\ell$;\\
  $(\mathcal{H}_{\ell,2})$\,\,\,
  $\lim\limits_{k\ra\infty}\varphi_{k}(s_k^{(j)})/\varphi_{k}(Lr_k^{(j)})=0$
  for all $1\leq j\leq\ell$ and all $L>0$;\\
  $(\mathcal{H}_{\ell,3})$\,\,\, $\lim\limits_{k\ra\infty}\Lambda_k(s_k^{(j)})=4\pi
  j$ for all $1\leq j\leq\ell$;\\
  $(\mathcal{H}_{\ell,4})$\,\,\,$\lim\limits_{L\ra\infty}\lim\limits_{k\ra\infty}
  \le(N_k(s_k^{(j-1)},r_k^{(j)}/L)+N_k(Lr_k^{(j)},s_k^{(j)})\ri)=0$ for all $1\leq
  j\leq\ell$.\\

  To prove Proposition 5.1, we follow the lines of
  \cite{L-R-S,Mar-Stru,Struwe}. Precisely we use induction
  as follows: $(\mathcal{H}_1)$ holds; if $(\mathcal{H}_\ell)$ holds,
  then either $(\mathcal{H}_{\ell+1})$ holds, or
  \be\label{0-}\lim_{L\ra\infty}\lim_{k\ra\infty}N_k(s_k^{(\ell)},\delta/L)=0.\ee
  In view of (\ref{energy}), we have
  \bna
  \Lambda_k(s_k^{(\ell)})&=&\int_{\mathbb{B}_{s_k^{(\ell)}}(\widetilde{x}_{i,k})}\widetilde{u}_k(x)
  \widetilde{f}_k(x,\widetilde{u}_k(x))dx\\
  &=&(1+o(1))\int_{\mathbb{B}_{s_k^{(\ell)}}(\widetilde{x}_{i,k})}\widetilde{u}_k(x)
  \widetilde{f}_k(x,\widetilde{u}_k(x))e^{\psi_i(x)}dx\\
  &=&(1+o(1))\int_{\phi_i^{-1}(\mathbb{B}_{s_k^{(\ell)}}(\widetilde{x}_{i,k}))}
  u_kf_k(x,u_k)dv_g\\
  &\leq& (1+o(1))\int_\Sigma u_kf_k(x,u_k)dv_g.
  \ena
  This together with (\ref{bdd}) and $(\mathcal{H}_{\ell,3})$ implies that the induction
  terminates after finitely-many  steps. Letting $\ell_0$ be the largest
  integer such that $(\mathcal{H}_{\ell_0})$ holds. Since $\widetilde{x}_{i,k}\ra 0$ as
  $k\ra\infty$,
  in view of the last assertion of Proposition 3.1,
  for any fixed $L>2/\delta$,
  \be\label{away}\lim_{k\ra\infty}\le\|\widetilde{u}_k-\widetilde{u}_\infty\ri\|
  _{C^1(\Omega\setminus\mathbb{B}_{\delta/L}(\widetilde{x}_{i,k}))}=0.\ee
     Moreover it follows from
  $(\mathcal{H}_{\ell_0,3})$ and (\ref{0-}) (with $\ell$ replaced by $\ell_0$) that
  \be\label{inner}
  \lim_{L\ra\infty}\lim_{k\ra\infty}\int_{\phi_i^{-1}(\mathbb{B}_{\delta/L}(\widetilde{x}_{i,k}))}u_kf_k(x,u_k)dv_g=
  \lim_{L\ra\infty}\lim_{k\ra\infty}\int_{\mathbb{B}_{\delta/L}(\widetilde{x}_{i,k})}\widetilde{u}_k\widetilde{f}_k(x,\widetilde{u}_k)dx
  =4\pi\ell_0.
  \ee
  Recalling equation (\ref{e1}), we obtain (\ref{identity}) by combining
  (\ref{away}) and (\ref{inner}) with $I^{(i)}=\ell_0$, and thus complete the proof of Proposition 5.1.\\

  The proof of the above induction process will be divided into the two steps below.\\

  {\it Step 1. The property $(\mathcal{H}_1)$ holds.}\\

  For any function $h:\Omega=\phi_i(U_i)\ra\mathbb{R}$, denote the
  spherical average of $h$ around $\widetilde{x}_{i,k}$ by
  $$\overline{h}(r)=\f{1}{2\pi r}\int_{\p\mathbb{B}_r(\widetilde{x}_{i,k})}hd\sigma,\quad\forall 0<r<\delta.$$
  Let $w_k$  be the unscaled function with respect to the blow-up sequence $\eta_{i,k}$ as in
  (\ref{blow-up-funct}), namely
  $$w_k(x)=u_k(x_{i,k})(\widetilde{u}_k(x)-u_k(x_{i,k})),\quad x\in\Omega.$$
  The decay estimate on $\overline{w}_k$ near the point $\widetilde{x}_{i,k}$ is
  crucial
  for the property $(\mathcal{H}_1)$. Precisely we have the following result.\\

  \noindent{\bf Lemma 5.2} {\it Given $0<\epsilon<1$. Let $T_k$
  be the smallest number such that
  $\varphi_{k}(T_k)=\epsilon u_k(x_{i,k})$. Then
  ${r_{i,k}}/{T_k}\ra 0$ as $k\ra\infty$, where $r_{i,k}$ is as in (\ref{scal}).
  Moreover, for any $b<2$, there exist some integer $k_0$ and a constant $C$
  such that when $k\geq k_0$, we have
  \be\label{w}\overline{w}_k(r)\leq b\log\f{r_{i,k}}{r}+C\ee
  for all $0\leq r\leq T_k$ and
  \be\label{n}\lim_{k\ra\infty}\Lambda_k(T_k)=4\pi.\ee}
  {\it Proof.} It follows from Proposition 3.1 and the definition of $T_k$
  that ${r_{i,k}}=o({T_k})$ as $k\ra\infty$. In view of (\ref{e1}), $\widetilde{u}_k$
  satisfies the equation
  \be\label{u-tild}-\Delta_{\mathbb{R}^2}\widetilde{u}_k=e^{\psi_i}(\widetilde{f}_k(x,\widetilde{u}_k)-\widetilde{\tau}_k
  \widetilde{u}_k)\quad{\rm in}\quad \Omega.\ee
  Let $(v_k)$ be a sequence of
  solutions
  to
  \be\label{nh}\le\{\begin{array}{lll}-\Delta_{\mathbb{R}^2}v_k=e^{\psi_i} \widetilde{f}_k(x,\widetilde{u}_k)
  &{\rm on} &\mathbb{B}_{T_k}(\widetilde{x}_{i,k})\\
  [1.5ex] v_k=\widetilde{u}_k&{\rm
  on}&\p\mathbb{B}_{T_k}(\widetilde{x}_{i,k}).\end{array}\ri.\ee
  Then we have by (\ref{u-tild})
  \be\label{nhh}\le\{\begin{array}{lll}-\Delta_{\mathbb{R}^2}(v_k-\widetilde{u}_k)=e^{\psi_i}
  \widetilde{\tau}_k\widetilde{u}_k
  &{\rm on} &\mathbb{B}_{T_k}(\widetilde{x}_{i,k})\\
  [1.5ex] v_k-\widetilde{u}_k=0&{\rm
  on}&\p\mathbb{B}_{T_k}(\widetilde{x}_{i,k}).\end{array}\ri.\ee
  Applying elliptic estimates to (\ref{nhh}), we can find some
  constant $C$ independent of $k$ such
  that
  $$|v_k(x)-\widetilde{u}_k(x)|\leq C\,\,\,{\rm for\,\,\,all}\,\,\,
  x\in \mathbb{B}_{T_k}(\widetilde{x}_{i,k}).$$
  Moreover, it follows from Proposition 4.1 that
  \be\label{inf}\inf_{\p\mathbb{B}_{T_k}(\widetilde{x}_{i,k})}\widetilde{u}_k\geq
  \varphi_{k}(T_k)-C\ee
  for some constant $C$ depending only on the Riemannian metric $g$. Applying the
  maximum principle to (\ref{nh}), we have by (\ref{inf})
  \be\label{lower}\widetilde{u}_k(x)\geq \varphi_{k}(T_k)-C\,\,\,{\rm for\,\,\,all}\,\,\,
  x\in \mathbb{B}_{T_k}(\widetilde{x}_{i,k}).\ee
  Note that $\varphi_{k}(T_k)=\epsilon u_k(x_{i,k})$. For any $0\leq t\leq
  T_k$, we have by (\ref{lower}) and the fact that $u_k\ra u_\infty$ strongly in
  $L^2(\Sigma)$
  \be\label{err} u_k(x_{i,k})\int_{\mathbb{B}_t(\widetilde{x}_{i,k})} e^{\psi_i}\widetilde{\tau}_k\widetilde{u}_kdx
  \leq
  \f{\|\tau_k\|_{L^\infty(\Sigma)}}{\epsilon}\int_{\mathbb{B}_t(\widetilde{x}_{i,k})}e^{\psi_i}
  \le(\widetilde{u}_k^2+C\widetilde{u}_k\ri)dx=o(1).\ee
  For any $Lr_{i,k}\leq t\leq T_k$, we obtain by Proposition 3.1
  \bea\label{err1}{\nonumber}
  -u_k(x_{i,k})\int_{\mathbb{B}_{t}(\widetilde{x}_{i,k})}e^{\psi_i}\widetilde{f}_k(x,\widetilde{u}_k)dx
  &\leq&-u_k(x_{i,k})\int_{\mathbb{B}_{Lr_{i,k}}(\widetilde{x}_{i,k})}e^{\psi_i}\widetilde{f}_k(x,\widetilde{u}_k)dx\\
  &=&-r_{i,k}^{-2}\int_{\mathbb{B}_{Lr_{i,k}}(\widetilde{x}_{i,k})}
  e^{\psi_i}\f{\widetilde{f}_k(x,\widetilde{u}_k)}{f_k(x_{i,k},u_k(x_{i,k}))}dx{\nonumber}\\{\nonumber}
  &=&-(1+o(1))\int_{\mathbb{B}_L(0)}e^{(2+o(1))\eta_\infty}dx\\
  &=&-4\pi+o(1),\eea
  where $o(1)\ra 0$ as $k\ra\infty$ first, and then $L\ra\infty$.
  In view of (\ref{u-tild}), $\overline{w}_k$ satisfies
  $$-\Delta_{\mathbb{R}^2}\overline{w}_k= u_k(x_{i,k})\overline{e^{\psi_i}\widetilde{f}_k(x,\widetilde{u}_k)}
  -u_k(x_{i,k}) \overline{e^{\psi_i} \widetilde{\tau}_k\widetilde{u}_k}.$$
  Then we have for any $Lr_{i,k}\leq t\leq T_k$
  \bna
  2\pi
  t\overline{w}_k^\prime(t)&=&\int_{\p\mathbb{B}_{t}(\widetilde{x}_{i,k})}\p_\nu\overline{w}_k
  d\sigma=\int_{\mathbb{B}_{t}(\widetilde{x}_{i,k})}\Delta_{\mathbb{R}^2}\overline{w}_k
  dx\\
  &=&-u_k(x_{i,k})\int_{\mathbb{B}_{t}(\widetilde{x}_{i,k})}\overline{e^{\psi_i}\widetilde{f}_k(x,\widetilde{u}_k)}dx
  +u_k(x_{i,k})\int_{\mathbb{B}_{t}(\widetilde{x}_{i,k})}\overline{e^{\psi_i}\widetilde{\tau}_k
  \widetilde{u}_k}dx\\
  &=&-u_k(x_{i,k})\int_{\mathbb{B}_{t}(\widetilde{x}_{i,k})}{e^{\psi_i}\widetilde{f}_k(x,\widetilde{u}_k)}dx
  +u_k(x_{i,k})\int_{\mathbb{B}_{t}(\widetilde{x}_{i,k})}{e^{\psi_i}\widetilde{\tau}_k
  \widetilde{u}_k}dx\\
  &\leq& -4\pi+o(1).
  \ena
  Here we used (\ref{err}) and (\ref{err1}) in the last inequality. Thus for any
  $b<2$, there exists some integer $k_0$ such that
  $$\overline{w}_k^\prime(t)\leq -\f{b}{t}\quad{\rm for\,\,all}\quad k\geq k_0.$$
  This together with Proposition 3.1 leads to
  \bna
  \overline{w}_k(t)&\leq&\overline{w}_k(Lr_{i,k})-b\log\f{t}{Lr_{i,k}}{\nonumber}\\
  &\leq&\log\f{1}{1+L^2}-b\log\f{t}{Lr_{i,k}}+o(1){\nonumber}\\
  &\leq& b\log\f{r_{i,k}}{t}+C
  \ena
  for some constant $C$, all $k\geq k_0$, and all $Lr_{i,j}\leq t\leq
  T_k$. It follows from Proposition 3.1 again that the above inequality
  also holds for $0\leq t\leq Lr_{i,k}$. Hence (\ref{w}) holds.

  By (\ref{w}) and (\ref{lower}) we have
  \bna
  (\epsilon-1)u_k^2(x_{i,k})-Cu_k(x_{i,k})\leq\overline{w}_k(r)
  \leq C,\quad \forall r\in
  [Lr_{i,k}, T_k].
  \ena
  Hence there holds for $Lr_{i,k}\leq r\leq T_k$
  \bea\nonumber
  \varphi_{k}^2(r)-u_k^2(x_{i,k})&=&\le(1+\f{\varphi_{k}(r)}{u_k(x_{i,k})}\ri)\overline{w}_k(r)\\
  \nonumber&=&\le(2+\f{\overline{w}_k}{u_k^2(x_{i,k})}\ri)\overline{w}_k(r)\\
  \nonumber&\leq&(1+\epsilon+o(1))\overline{w}_k(r)+(1-\epsilon+o(1))C\\
  \label{phi-r}&\leq& (1+2\epsilon/3)b\log\f{r_{i,k}}{r}+C,
  \eea
  provided that $k$ is sufficiently large.
  For $0<r<\delta$ we denote
  \be\label{theta}\theta_k(r)=\theta_k^{(i)}(r)=\f{1}{2\pi r}\int_{\p\mathbb{B}_r(\widetilde{x}_{i,k})}
  \widetilde{f}_k(x,\varphi_k(r))d\sigma.\ee
  Taking $b$ such that
  $(1+2\epsilon/3)b=2+\epsilon$ in (\ref{phi-r}) and recalling (H4) and (H5), we
  can find some constant $C$ such that
  for $Lr_{i,k}\leq r\leq T_k$
  \bea\nonumber
  \f{\theta_{k}(r)}{f_k(x_{i,k},u_k(x_{i,k}))}&=&\f{\theta_{k}(r)}
  {\widetilde{f}_k(\widetilde{x}_{i,k},\varphi_{k}(r))}
  \f{\widetilde{f}_k(\widetilde{x}_{i,k},\varphi_{k}(r))}{f_k(x_{i,k},u_k(x_{i,k}))}\\
  \nonumber&=&(1+o(1))\f{{f}_k({x}_{i,k},\varphi_{k}(r))}{f_k(x_{i,k},u_k(x_{i,k}))}\\
  \nonumber&=&(1+o(1))e^{(1+o(1))(\varphi_{k}^2(r)-u_k^2(x_{i,k}))}\\
  \label{thk}&\leq&C\le(\f{r_{i,k}}{r}\ri)^{2+\epsilon}
  \eea
  for sufficiently large $k$. For $0<s<t<\delta$, we define next a function analogous to  (\ref{neck}) as below.
  \be\label{neck-symm}\overline{N}_k(s,t)=\overline{N}_k^{(i)}(s,t)=2\pi\int_s^tr\varphi_k(r)\theta_k(r)dr.\ee
   In view of (\ref{w}) and (\ref{thk}), we estimate
  \bna
  \overline{N}_k(Lr_{i,k},T_k)&=&2\pi\int_{Lr_{i,k}}^{T_k}r\varphi_k(r)\theta_k(r)dr\\
  &=&2\pi r_{i,k}^{-2}\int_{Lr_{i,k}}^{T_k}r\f{\varphi_k(r)}{u_k(x_{i,k})}\f{\theta_k(r)}
  {f_k(x_{i,k},u_k(x_{i,k}))}dr\\
  &\leq&2\pi(1+o(1))Cr_{i,k}^{\epsilon}\int_{Lr_{i,k}}^{T_k}\f{1}{r^{1+\epsilon}}dr\\
  &\leq&2\pi(1+o(1))C\epsilon^{-1}L^{-\epsilon}.
  \ena
  This leads to
  \be\label{nn}\lim_{L\ra\infty}\lim_{k\ra\infty}\overline{N}_k(Lr_{i,k},T_k)=0.\ee
  Since Proposition 4.1 implies that
  $$u_k^2(x)-\varphi_{k}^2(r)\leq C\quad{\rm for \,\,all}\quad x\in\p\mathbb{B}_{r}(\widetilde{x}_{i,k}),$$
  there holds
  $$N_k(Lr_{i,k},T_k)\leq C\overline{N}_k(Lr_{i,k},T_k)+o(1).$$
  This together with (\ref{nn}) leads to
  \be\label{m1}\lim_{L\ra\infty}\lim_{k\ra\infty}{N}_k(Lr_{i,k},T_k)=0.\ee
  By Proposition 3.1,
  \bna
  \Lambda_k(Lr_{i,k})=\int_{\mathbb{B}_{Lr_{i,k}}(\widetilde{x}_{i,k})}
  \widetilde{u}_k\widetilde{f}_k(x,\widetilde{u}_k)dx
  =(1+o(1))\int_{\mathbb{B}_L(0)}e^{2\eta_\infty}dx.
  \ena
  Hence
  \be\label{m2}\lim_{L\ra\infty}\lim_{k\ra\infty}\Lambda_k(Lr_{i,k})=4\pi.\ee
  Thus (\ref{n}) follows  immediately from (\ref{m1}) and (\ref{m2}).
  $\hfill\Box$\\

  By Lemma 5.2 we may choose a subsequence $u_k$, numbers $\epsilon_k\searrow
  0$ as $k\ra\infty$ and $s_k=T_k(\epsilon_k)$ with $r_{i,k}/s_k\ra
  0$, $\varphi_{k}(s_k)\ra\infty$ as $k\ra\infty$ and such that
  $$\lim_{k\ra\infty}\Lambda_k(s_k)=4\pi,\quad \lim_{L\ra\infty}\lim_{k\ra\infty}{N}_k(Lr_{i,k},s_k)=0,$$
  while in addition
  $$\lim_{k\ra\infty}\f{\varphi_{k}(s_k)}{\varphi_{k}(Lr_{i,k})}=0,\quad\forall\, L>0.$$
  Let $r_k^{(1)}=r_{i,k}$, $s_k^{(1)}=s_k$. Then $(\mathcal{H}_1)$
  holds and Step 1 is finished. \\

  \noindent{\it Step 2. Suppose that $(\mathcal{H}_\ell)$ already holds for some integer $\ell\geq 1$,
  namely there exist sequences $s_k^{(0)}=0<r_k^{(1)}<s_k^{(1)}<\cdots<r_k^{(\ell)}<s_k^{(\ell)}=o(1)$
  such that $(\mathcal{H}_{\ell,1})$ up to
  $(\mathcal{H}_{\ell,4})$ hold. Then we shall prove that either
  $\lim_{L\ra\infty}\lim_{k\ra\infty}N_k(s_k^{(\ell)},\delta/L)=0$
  or $(\mathcal{H}_{\ell+1})$ holds.}\\

  Setting
  \be\label{Pt}P_k(t)=P_k^{(i)}(t)=t\int_{\p\mathbb{B}_t}\widetilde{u}_k\widetilde{f}_k(x,\widetilde{u}_k)d\sigma,
  \quad \overline{P}_k(t)=\overline{P}_k^{(i)}(t)=2\pi t^2\varphi_k(t)\theta_k(t)\ee
  and assuming $(\mathcal{H}_\ell)$ holds,
  we have the following result.\\

  \noindent{\bf Lemma 5.3} {\it There exists a constant $C_0$ depending only on
  the upper bound of the total energy (\ref{bdd}) and the Riemannian metric $g$ such that for $s_k^{(\ell)}\leq t_k=o(1)$, there holds
  \be\label{lemm}\overline{N}_{k}(s_k^{(\ell)},t_k)\leq \overline{P}_{k}(t_k)+C_0
  \overline{N}_{k}^2(s_k^{(\ell)},t_k)+o(1),\ee
  where $o(1)\ra 0$ as $k\ra\infty$, $\overline{N}_{k}$ and $\overline{P}_{k}$ are defined as in
  (\ref{neck-symm}) and (\ref{Pt}) respectively.}\\

  \noindent{\it Proof.}  We first claim that
  there exists a constant $C$ depending only on $\delta$ and the Riemannian metric $g$
  such that
  \be\label{cla}\varphi_{k}(s)\leq \sup_{\p\mathbb{B}_s(\widetilde{x}_{i,k})}\widetilde{u}_k
  \leq \inf_{\p\mathbb{B}_r(\widetilde{x}_{i,k})}\widetilde{u}_k+C\leq \varphi_{k}(r)+C
  \,\,\,{\rm for\,\,all}\,\,\, 0<r<s\leq \delta.\ee
  To see the last inequality, we set $v_k$ be a positive solution of
  \be\label{nh2}\le\{\begin{array}{lll}-\Delta_{\mathbb{R}^2}v_k=e^{\psi_i} \widetilde{f}_k(x,\widetilde{u}_k)
  &{\rm in}&\mathbb{B}_{\delta}(\widetilde{x}_{i,k})\\
  [1.5ex] v_k=\widetilde{u}_k&{\rm
  on}&\p\mathbb{B}_{\delta}(\widetilde{x}_{i,k}).\end{array}\ri.\ee
  Thus we have by (\ref{e1})
  \be\label{nhh2}\le\{\begin{array}{lll}-\Delta_{\mathbb{R}^2}(v_k-\widetilde{u}_k)=e^{\psi_i}
  \widetilde{\tau}_k\widetilde{u}_k
  &{\rm in} &\mathbb{B}_{\delta}(\widetilde{x}_{i,k})\\
  [1.5ex] v_k-\widetilde{u}_k=0&{\rm
  on}&\p\mathbb{B}_{\delta}(\widetilde{x}_{i,k}).\end{array}\ri.\ee
  Noting that $\|e^{\psi_i}
  \widetilde{\tau}_k\widetilde{u}_k\|_{L^p(\mathbb{B}_{\delta}(\widetilde{x}_{i,k}))}$ is
  bounded for any $p>1$ and
  applying elliptic regularity estimates to (\ref{nhh2}), we then find some
  constant $C=C(\delta)$ such that
  \be\label{uk}
  v_k(x)-C\leq \widetilde{u}_k(x)\leq v_k(x)+C\quad{\rm
  for\,\,all}\quad x\in\mathbb{B}_{\delta}(\widetilde{x}_{i,k}).
  \ee
   By (\ref{nh2}), we have for $0<r<\delta$
    $$-(r\overline{v}_k^\prime(r))^\prime=r\,\overline{e^{\psi_i} \widetilde{f}_k(x,\widetilde{u}_k)}.$$
  Integration from $0$ to $r$ gives
  $$-r\overline{v}_k^\prime(r)=\int_0^r r\,\overline{e^{\psi_i} \widetilde{f}_k(x,\widetilde{u}_k)}dr.$$
  Hence
  \be\label{mono}\overline{v}_k^\prime(r)\leq 0\quad {\rm for\,\,all}\quad 0<r<\delta.\ee
  Now fix $0<r<s\leq\delta$. There exist two points
  $\xi\in\p\mathbb{B}_{r}(\widetilde{x}_{i,k})$ and
  $\zeta\in\p\mathbb{B}_{s}(\widetilde{x}_{i,k})$ such that
  $$v_k(\xi)=\overline{v}_k(r),\quad v_k(\zeta)=\overline{v}_k(s).$$
  This together with the gradient estimate (Proposition 4.1), (\ref{uk}),
  and (\ref{mono}) leads to
  \bna
  \sup_{\p\mathbb{B}_s(\widetilde{x}_{i,k})}\widetilde{u}_k&\leq& \widetilde{u}_k(\zeta)+C
  \leq v_k(\zeta)+C\\&\leq& v_k(\xi)+C
  \leq\inf_{\p\mathbb{B}_r(\widetilde{x}_{i,k})}\widetilde{u}_k+C.
  \ena
  This confirms our claim (\ref{cla}).

  Next we calculate
  \bna\theta_{k}^\prime(r)&=&\f{d}{dr}\le(\f{1}{2\pi}\int_0^{2\pi}\widetilde{f}_k
  \le(\widetilde{x}_{i,k}^1+r\cos\theta,\widetilde{x}_{i,k}^2+r\sin\theta,\varphi_{k}(r)\ri)d\theta\ri)\\
  &=&\f{1}{2\pi}\int_0^{2\pi}\nabla_x\widetilde{f}_k
  \le(\widetilde{x}_{i,k}^1+r\cos\theta,\widetilde{x}_{i,k}^2+r\sin\theta,\varphi_{k}(r)\ri)\cdot
  (\cos\theta,\sin\theta)d\theta\\
  &&+\f{1}{2\pi}\int_0^{2\pi}\widetilde{f}_k^\prime
  \le(\widetilde{x}_{i,k}^1+r\cos\theta,\widetilde{x}_{i,k}^2+r\sin\theta,\varphi_{k}(r)\ri)
  \varphi_k^\prime(r)d\theta,\ena
  where we write
  $\widetilde{x}_{i,k}=(\widetilde{x}_{i,k}^1,\widetilde{x}_{i,k}^2)$.
  In view of (H4), we obtain
  \be\label{der}|\theta_{k}^\prime(r)|\leq C\le(1+\theta_{k}(r)+\varphi_{k}(r)
  |\varphi_{k}^\prime(r)|\theta_{k}(r)\ri).\ee
  For $s=s_k^{(\ell)}\leq t\leq t_k$, we have by  equation (\ref{e1}).
  \bna
  -2\pi
  t\varphi_{k}^\prime(t)&=&-\int_{\p\mathbb{B}_{t}(\widetilde{x}_{i,k})}\p_\nu\varphi_{k}
  d\sigma=-\int_{\mathbb{B}_{t}(\widetilde{x}_{i,k})}\Delta_{\mathbb{R}^2}\varphi_{k}dx\\
  &=&\int_{\mathbb{B}_{t}(\widetilde{x}_{i,k})}\overline{e^{\psi_i}\widetilde{f}_k(x,\widetilde{u}_k)}dx
  -\int_{\mathbb{B}_{t}(\widetilde{x}_{i,k})}\overline{e^{\psi_i}\widetilde{\tau}_k\widetilde{u}_k}dx\\
  &=&\int_{\mathbb{B}_{t}(\widetilde{x}_{i,k})}e^{\psi_i}\widetilde{f}_k(x,\widetilde{u}_k)dx
  -\int_{\mathbb{B}_{t}(\widetilde{x}_{i,k})}e^{\psi_i}\widetilde{\tau}_k\widetilde{u}_kdx.
  \ena
  It follows from (H4), (H5) and Proposition 4.1 that
  $\widetilde{f}_k(x,\widetilde{u}_k)\leq
  C(1+\widetilde{f}_k(x,\varphi_k(r)))\leq C(1+\theta_k(r))$,
  where $r=|x-\widetilde{x}_{i,k}|$. Combining (H4), (H5) and
  (\ref{cla}), we have
  $$\int_{\mathbb{B}_s(\widetilde{x}_{k})}\varphi_k(s)\theta_k(r)dx\leq C(1+\Lambda_k(s)),$$
  where we used $r=|x-\widetilde{x}_{i,k}|$.
  Note that $\varphi_k(s)\ra\infty$ as $k\ra\infty$. We then obtain
  \bna
  -2\pi
  t\varphi_{k}^\prime(t)&\leq&
  C\int_{\mathbb{B}_{t}(\widetilde{x}_{i,k})}\le(1+\theta_{k}(r)\ri)dx\\
  &\leq&C\int_{\mathbb{B}_{t}(\widetilde{x}_{i,k})\setminus\mathbb{B}_{s}
  (\widetilde{x}_{i,k})}\theta_{k}(r)dx+\f{C}{\varphi_{k}(s)}
  \int_{\mathbb{B}_{s}(\widetilde{x}_{i,k})}\varphi_{k}(s)\theta_{k}(r)dx+o(1)\\
    &\leq&C\overline{N}_k(s,t)+o(1).
  \ena
  This immediately leads to
  \be\label{dr1}
  -\pi\int_s^tr^2\varphi_{k}^\prime(r)\theta_{k}(r)dr\leq
  C\overline{N}_k^2(s,t)+o(1).
  \ee
  Similarly we have
  \bea{\nonumber}
  -2\pi
  t\varphi_{k}(t)\varphi_{k}^\prime(t)&=&-\int_{\p\mathbb{B}_{t}(\widetilde{x}_{i,k})}\varphi_{k}(t)\p_\nu\varphi_{k}
  d\sigma=-\int_{\mathbb{B}_{t}(\widetilde{x}_{i,k})}\varphi_{k}(t)\Delta_{\mathbb{R}^2}\varphi_{k}dx\\
  {\nonumber}&=&\int_{\mathbb{B}_{t}(\widetilde{x}_{i,k})}\varphi_{k}(t)\overline{e^{\psi_i}\widetilde{f}_k(x,\widetilde{u}_k)}dx
  -\int_{\mathbb{B}_{t}(\widetilde{x}_{i,k})}\varphi_{k}(t)\overline{e^{\psi_i}\widetilde{\tau}_k\widetilde{u}_k}dx\\
  \label{tph}&=&\int_{\mathbb{B}_{t}(\widetilde{x}_{i,k})}\varphi_{k}(t)\overline{e^{\psi_i}\widetilde{f}_k(x,\widetilde{u}_k)}dx+o(1),
  \eea
  where the last equality follows from (\ref{cla}) and $u_k\ra u_\infty$ strongly in
  $L^2(\Sigma)$. Repeatedly using (\ref{cla}), we obtain
  \bna
  \int_{\mathbb{B}_{t}(\widetilde{x}_{i,k})}\varphi_{k}(t)\overline{e^{\psi_i}\widetilde{f}_k(x,\widetilde{u}_k)}dx
  &\leq&C\int_{\mathbb{B}_{t}(\widetilde{x}_{i,k})}\varphi_{k}(t)\le(1+\theta_{k}(r)\ri)dx\\
  &\leq&C\int_{\mathbb{B}_{t}(\widetilde{x}_{i,k})\setminus\mathbb{B}_{s}(\widetilde{x}_{i,k})}\le(1+\varphi_{k}(r)\ri)
  \le(1+\theta_{k}(r)\ri)dx\\
  &&+C\int_{\mathbb{B}_{s}(\widetilde{x}_{i,k})\setminus\mathbb{B}_{Lr_k^{(\ell)}}(\widetilde{x}_{i,k})}\le(1+\varphi_{k}(r)\ri)
  \le(1+\theta_{k}(r)\ri)dx\\
  &&+C\int_{\mathbb{B}_{Lr_k^{(\ell)}}(\widetilde{x}_{i,k})}\le(1+\varphi_{k}(s)\ri)
  \le(1+\theta_{k}(r)\ri)dx\\
  &\leq&C\le(\overline{N}_k(s,t)+\overline{N}_k(Lr_k^{(\ell)},s)+\f{\varphi_{k}(s)}{\varphi_{k}(Lr_k^{(\ell)})}
  \le({\Lambda}_k(Lr_k^{(\ell)})+o(1)\ri)\ri).
  \ena
  This together with
  (\ref{tph}), $(\mathcal{H}_{\ell,2})$ and $(\mathcal{H}_{\ell,4})$ implies
  \be\label{t-p}2\pi t\varphi_{k}(t)|\varphi_{k}^\prime(t)|\leq
  C\overline{N}_k(s,t)+o(1).\ee
  Obviously
  $$\int_s^tr^2\varphi_{k}(r)dr=o(1),\quad \int_s^tr^2\varphi_{k}(r)\theta_{k}(r)dr=o(1).$$
  It then follows from (\ref{der}) and (\ref{t-p}) that
  \bea\nonumber
  -\pi\int_s^tr^2\varphi_{k}(r)\theta_{k}^\prime(r)dr&\leq&\pi C\int_s^tr^2\varphi_{k}^2(r)
  |\varphi_{k}^\prime(r)|\theta_{k}(r)dr+o(1)\\
  &\leq&C\overline{N}_k^2(s,t)+o(1).\label{i-2}
  \eea
  Integration by parts gives
  \bna
  \overline{N}_k(s,t)&=&\int_s^t2\pi
  r\varphi_{k}(r)\theta_{k}(r)dr\\
  &\leq&\pi t^2\varphi_{k}(t)\theta_{k}(t)-\pi\int_s^tr^2\varphi_{k}^\prime(r)\theta_{k}(r)dr
  -\pi\int_s^tr^2\varphi_{k}(r)\theta_{k}^\prime(r)dr.
  \ena
  This together with (\ref{dr1}) and (\ref{i-2}) implies (\ref{lemm}). $\hfill\Box$\\

  \noindent{\bf Lemma 5.4} {\it Let $C_0$ be the constant as in Lemma 5.3. Let $t_k$ be such that for a subsequence
  $$s_k^{(\ell)}<t_k=o(1),\quad 0<\lim_{k\ra\infty}\overline{N}_k(s_k^{(\ell)},t_k)=\alpha<\f{1}{2C_0}.$$
  Then $s_k^{(\ell)}=o(t_k)$ as $k\ra\infty$, $\liminf\limits_{k\ra\infty}\overline{P}_k(t_k)\geq \alpha/2$, and
  \be\label{3-}\lim_{L\ra\infty}\lim_{k\ra\infty}\overline{N}_k(s_k^{(\ell)},t_k/L)=0,\ee
  where $\overline{N}_k$ and $\overline{P}_k$ are as defined in (\ref{neck-symm}) and (\ref{Pt})
  respectively.}\\

  \noindent{\it Proof.} We first claim that
  \be\label{1-}\lim_{L\ra\infty}\lim_{k\ra\infty}\overline{N}_k(s_k^{(\ell)},Ls_k^{(\ell)})=0.\ee
  Actually, in view of (\ref{cla}), we have for $0<t\leq t_k$
  \be\label{P-N}\overline{P}_k(t)\leq C\overline{N}_k(t/2,t)+o(1)\leq C\overline{P}_k(t/2)+o(1),\ee
  and
  $$\overline{N}_k(t,2t)\leq C\overline{N}_k(t/2,t)+o(1).$$
  In particular, for any $j\in\mathbb{N}$ there holds
  \bna
  \lim_{k\ra\infty}\overline{N}_k(2^{j-1}s_k^{(\ell)},2^{j}s_k^{(\ell)})
  &\leq&C\lim_{k\ra\infty}\overline{N}_k(2^{j-2}s_k^{(\ell)},2^{j-1}s_k^{(\ell)})\\
  &\leq&C^{j}\lim_{k\ra\infty}\overline{N}_k(s_k^{(\ell)}/2,s_k^{(\ell)})=0.
  \ena
  If $L\leq 2^j$, we obtain
  $$\lim_{k\ra\infty}\overline{N}_k(s_k^{(\ell)},Ls_k^{(\ell)})\leq \lim_{k\ra\infty}\sum_{m=1}^j
  \overline{N}_k(2^{m-1}s_k^{(\ell)},2^ms_k^{(\ell)})=0.$$
  Thus our claim (\ref{1-}) follows immediately. One can see from (\ref{1-}) that $s_k^{(\ell)}/t_k\ra
  0$ as $k\ra\infty$. By Lemma 5.3,
  \be\label{2-}
  \liminf_{k\ra\infty}\overline{P}_k(t_k)\geq
  \f{1}{2}\lim_{k\ra\infty}\overline{N}_k(s_k^{(\ell)},t_k)=\f{\alpha}{2}.
  \ee
  Now we show (\ref{3-}).
  Assuming the contrary, there holds
  $$\lim_{L\ra\infty}\lim_{k\ra\infty}\overline{N}_k(s_k^{(\ell)},t_k/L)=\beta>0.$$
  Then we have for any fixed $L\geq 1$ and all sufficiently large $k$
  $$\f{\beta}{2}\leq \overline{N}_k(s_k^{(\ell)},t_k/L)\leq\overline{N}_k(s_k^{(\ell)},t_k)<\f{1}{2C_0}.$$
  Applying (\ref{lemm}) with $t_k/L$ instead of $t_k$, we get
  $$\lim_{k\ra\infty}\overline{P}_k(t_k/L)\geq \f{\beta}{4},$$
  and then by (\ref{P-N})
  $$C\lim_{k\ra\infty}\overline{N}_k\le(t_k/(2L),t_k/L\ri)\geq\lim_{k\ra\infty}\overline{P}_k(t_k/L)\geq \f{\beta}{4}.$$
  Choosing $L=2^m$, $m=0,1,\cdots,j-1$, we have
  \bna
  \f{j\beta}{4}\leq C\lim_{k\ra\infty}\overline{N}_k(2^{-j}t_k,t_k)
  \leq C(1+\limsup_{k\ra\infty}{\Lambda}_k(t_k))\leq C.
  \ena
  We get a contradiction by letting $j\ra\infty$ and obtain
  (\ref{3-}). $\hfill\Box$\\

  \noindent{\bf Lemma 5.5} {\it Suppose that
  \be\label{caseA}\lim_{k\ra\infty}\sup_{s_k^{(\ell)}<t<t_k}\overline{P}_k(t)=0
  \quad{\rm for\,\,any\,\,sequence}\,\, t_k\ra 0\,\,{\rm as}\,\, k\ra \infty.\ee
  Then we have
  $$\lim_{L\ra\infty}\lim_{k\ra\infty}\overline{N}_k(s_k^{(\ell)},\delta/L)=0.$$
  }
  {\it Proof.} In view of Lemma 5.4, it suffices to prove
  \be\label{cd}\lim_{L\ra\infty}\lim_{k\ra\infty}\sup_{s_k^{(\ell)}<t<\delta/L}\overline{P}_k(t)=0.\ee
  Indeed, if we take some number $t_{k,L}\in (s_k^{(\ell)},\delta/L)$ such
  that
  $$\overline{P}_k(t_{k,L})=\sup_{s_k^{(\ell)}<t<\delta/L}\overline{P}_k(t),$$
  then either
  \be\label{ze}\lim_{k\ra\infty}t_{k,L}=0,\ee
  or
  \be\label{ze1}\lim_{k\ra\infty}t_{k,L}=t_L^\ast>0.\ee
  In case of (\ref{ze}), we already have (\ref{cd}) because of
  (\ref{caseA}). While in case of (\ref{ze1}), we have by using (\ref{cla})
  \bea\label{eli}
  \overline{P}_k(t_{k,L})
  \leq Ct_{k,L}^2\le(1+\varphi_{k}(t_{L}^\ast/2)\theta_{k}(t_{L}^\ast/2)\ri)
   \eea
    for sufficiently large $k$.
  Note that
  $\p\mathbb{B}_{t_L^\ast/2}(\widetilde{x}_{i,k})\subset\mathbb{B}_{t_L^\ast}(\widetilde{x}_{i}^\ast)
  \setminus\mathbb{B}_{t_L^\ast/3}(\widetilde{x}_{i}^\ast)$
  for sufficiently large $k$, and that $t_{k,L}\leq \delta/L\ra 0$ as $k\ra\infty$ first and then $L\ra\infty$.
  Moreover, by Proposition 3.1, we have $u_k\ra u_\infty$ in
  $C^1_{\rm loc}(\Sigma\setminus \cup_{j=1}^N\{x_j^\ast\},\mathbb{R})$
  and $u_\infty\in C^1(\Sigma,\mathbb{R})$, In particular, $u_\infty$ is bounded on $\mathbb{B}_\delta(x_i^\ast)$.
  It then follows from (\ref{eli})
  that
  $$
  \lim_{L\ra\infty}\lim_{k\ra\infty}\overline{P}_k(t_{k,L})=0.
  $$
  Thus (\ref{cd}) holds again.  $\hfill\Box$\\

  If the assumption (\ref{caseA}) is not satisfied, then (\ref{P-N})
  implies that there exists a sequence $t_k\ra 0$ as $k\ra\infty$
  such that
  \be\label{(5)}\lim_{k\ra\infty}\overline{N}_k(s_k^{(\ell)},t_k)>0.\ee
  We shall show that the property $(\mathcal{H}_{\ell+1})$ holds.
  Take $r_k^{(\ell+1)}\in(s_k^{(\ell)},t_k)$ such that up to a
  subsequence, there holds
  $$0<\lim_{k\ra\infty}\overline{N}_k(s_k^{(\ell)},r_k^{(\ell+1)})<\f{1}{2C_0},$$
  where $C_0$ is as in Lemma 5.3. It then follows from Lemma 5.4
  that
  \bea
  \label{(6)}&&\lim_{k\ra\infty}\f{s_k^{(\ell)}}{r_k^{(\ell+1)}}=0,\quad
  \lim_{k\ra\infty}\overline{N}_k(s_k^{(\ell)},r_k^{(\ell+1)})>0,\\
  \label{(7)}&&\liminf_{k\ra\infty}\overline{P}_k(r_k^{(\ell+1)})>0,\quad
  \lim_{k\ra\infty}\varphi_{k}(r_k^{(\ell+1)})=\infty,
  \eea
  and that
  \be\label{(8)}\lim_{L\ra\infty}\lim_{k\ra\infty}\overline{N}_k(s_k^{(\ell)},r_k^{(\ell+1)}/L)=0.\ee
  Moreover, we have the following result.\\

  \noindent{\bf Lemma 5.6} {\it Up to a subsequence there holds
  $$\eta_k^{(\ell+1)}(x):=\varphi_{k}(r_k^{(\ell+1)})
  \le(\widetilde{u}_k(\widetilde{x}_{i,k}+r_k^{(\ell+1)} x)-\varphi_{k}(r_k^{(\ell+1)})\ri)
  \ra\eta^{(\ell+1)}(x)$$ in $C^1_{\rm loc}(\mathbb{R}^2\setminus\{0\})$ as $k\ra\infty$, where
  $$\eta^{(\ell+1)}(x)=\log\f{2}{\sqrt{\alpha_0}(1+|x|^2)}$$ and
  $$\int_{\mathbb{R}^2}e^{2\eta^{(\ell+1)}}dx=\f
{4\pi}{\alpha_0}$$
for some constant $\alpha_0>0$.}\\

  \noindent{\it Proof.} To simplify the notations we write $r_k=r_k^{(\ell+1)}$, $\eta_k=\eta_k^{(\ell+1)}$,
  and $\eta=\eta^{(\ell+1)}$.  For any fixed $L>0$, we set
  \be\label{(9)}v_k(x)=\widetilde{u}_k(\widetilde{x}_{i,k}+r_k x),
  \quad x\in\mathbb{B}_L(0)\setminus\mathbb{B}_{1/L}(0).\ee
  In view of Proposition 4.1, there exists some constant $C=C(L)$
  such that
  $$|\widetilde{u}_k^2(\widetilde{x}_{i,k}+r_k x)-\varphi_{k}^2(r_k)|\leq C,$$
  and thus
  \be\label{(11)}|\varphi_{k}(r_k)\le(\widetilde{u}_k(\widetilde{x}_{i,k}+r_k x)-\varphi_{k}(r_k)\ri)|\leq C.\ee
  Hence
  \be\label{(10)}\eta_k\quad{\rm is\,\,bounded\,\,in}\quad L^\infty_{\rm loc}(\mathbb{R}^2\setminus\{0\}).\ee
  Combining (\ref{(7)}) and (\ref{(11)}), we have
  $$v_k-\varphi_{k}(r_k)\ra 0\quad{\rm in}\quad L^\infty_{\rm loc}(\mathbb{R}^2\setminus\{0\})
  \quad{\rm as}\quad k\ra\infty,$$
  in particular
  \be\label{(13)}\f{v_k}{\varphi_{k}(r_k)}\ra 1\quad{\rm in}\quad L^\infty_{\rm loc}(\mathbb{R}^2\setminus\{0\})
  \quad{\rm as}\quad k\ra\infty.\ee
  By the equation (\ref{e1}), we write for $x\in\Omega_k=\{x\in\mathbb{R}^2:\widetilde{x}_{i,k}+r_k x
  \in \mathbb{B}_\delta(0)\}$
  \be\label{(14)}-\Delta_{\mathbb{R}^2}\eta_k(x)=e^{\psi_i(\widetilde{x}_{i,k}+r_k x)}\varphi_{k}(r_k)r_k^2
  \widetilde{f}(\widetilde{x}_{i,k}+r_k x,v_k(x))-e^{\psi_i(\widetilde{x}_{i,k}+r_k
  x)}\varphi_{k}(r_k)r_k^2\widetilde{\tau}_k(\widetilde{x}_{i,k}+r_kx)
  v_k(x).\ee
  Since $u_k\ra u_\infty$ strongly in $L^2(\Sigma)$, we have by
  using (\ref{(13)})
  \bea\nonumber
  r_k^2\varphi_{k}^2(r_k)&=&\f{r_k^2}{3\pi}\int_{\mathbb{B}_2(0)\setminus\mathbb{B}_{1}(0)}\varphi_{k}^2(r_k)dx\\
  \nonumber&=&(1+o(1))\f{r_k^2}{3\pi}\int_{\mathbb{B}_2(0)\setminus\mathbb{B}_{1}(0)}v_k^2(x)dx\\
  \nonumber&=&\f{1+o(1)}{3\pi}\int_{\mathbb{B}_{2r_k}(\widetilde{x}_{i,k})\setminus\mathbb{B}_{r_k}
  (\widetilde{x}_{i,k})}\widetilde{u}_k^2(y)dy\\
  \label{(15)}&\ra& 0\quad{\rm as}\quad k\ra\infty.
  \eea
  By (\ref{(7)}) we may assume
  \be\label{(16)}r_k^2\varphi_{k}^2(r_k)\theta_{k}(r_k)\ra \alpha_0>0.\ee
  Moreover, by (H4) and (H5) we have
  \bea
  \f{\widetilde{f}_k(\widetilde{x}_{i,k}+r_k
  x,v_k(x))}{\theta_{k}(r_k)}&=&(1+o(1))\f{\widetilde{f}_k(\widetilde{x}_{i,k}+r_k
  x,v_k(x))}{\widetilde{f}_k(\widetilde{x}_{i,k}+r_k
  x,\varphi_{k}(r_k))}\nonumber\\
  \nonumber&=&(1+o(1))e^{(1+o(1))(v_k^2(x)-\varphi_{i,k}^2(r_k))}\\\label{(17)}
  &=&(1+o(1))e^{(2+o(1))\eta_k(x)}.
  \eea
  Applying elliptic estimates to (\ref{(14)}), we conclude from
  (\ref{(10)}), (\ref{(15)})-(\ref{(17)}) that
  \be\label{(18)}\eta_k\ra\eta\quad{\rm in}\quad C^1_{\rm loc}(\mathbb{R}^2\setminus\{0\})
  \quad {\rm as}\quad k\ra\infty,\ee
  where $\eta$ satisfies
  \be\label{(19)}-\Delta_{\mathbb{R}^2}\eta=\alpha_0 e^{2\eta}\quad{\rm in}\quad \mathbb{R}^2\setminus\{0\}.\ee
  For any $L>0$, (\ref{(18)}) together with (\ref{bdd}), (\ref{(13)}) and (\ref{(16)})
  leads to
  \bna
  \int_{\mathbb{B}_L(0)\setminus\mathbb{B}_{1/L}(0)}e^{2\eta}dx&=&\lim_{k\ra\infty}
  \int_{\mathbb{B}_L(0)\setminus\mathbb{B}_{1/L}(0)}e^{2\eta_k}dx\\
  &=&\lim_{k\ra\infty}\int_{\mathbb{B}_L(0)\setminus\mathbb{B}_{1/L}(0)}
  \f{v_k(x)\widetilde{f}_k(\widetilde{x}_{i,k}+r_kx,v_k(x))}{\varphi_{k}(r_k)\theta_{k}(r_k)}dx\\
  &=&\f{1}{\alpha_0}\lim_{k\ra\infty}\int_{\mathbb{B}_{Lr_k}(\widetilde{x}_{i,k})\setminus
  \mathbb{B}_{r_k/L}(\widetilde{x}_{i,k})}\widetilde{u}_k(y)\widetilde{f}_k(y,\widetilde{u}_k(y))dy\\
  &\leq&\f{C}{\alpha_0}.
  \ena
  Letting $L\ra\infty$, we have
  $$\int_{\mathbb{R}^2}e^{2\eta}dx<\infty.$$
  It follows from (\ref{cla}), $(\mathcal{H}_{\ell,2})$ and $(\mathcal{H}_{\ell,4})$ that
  \bna
  \int_{\mathbb{B}_{s_k^{(\ell)}}(\widetilde{x}_{i,k})}\varphi_{k}(r_k)\widetilde{f}_k
  (y,\widetilde{u}_k(y))dy&\leq&\int_{\mathbb{B}_{s_k^{(\ell)}}(\widetilde{x}_{i,k})\setminus
  \mathbb{B}_{Lr_k^{(\ell)}}(\widetilde{x}_{i,k})}\varphi_{k}(r)\widetilde{f}_k(y,\widetilde{u}_k(y))dy\\
  &&+\f{\varphi_{k}(s_k^{(\ell)})}{\varphi_{k}(Lr_k^{(\ell)})}\int_{\mathbb{B}_{Lr_k^{(\ell)}}(\widetilde{x}_{i,k})}
  \varphi_{k}(r)\widetilde{f}_k(y,\widetilde{u}_k(y))dy+o(1)\\
  &\leq&
  N_k(Lr_k^{(\ell)},s_k^{(\ell)})+\f{\varphi_{k}(s_k^{(\ell)})}{\varphi_{k}(Lr_k^{(\ell)})}\Lambda_k(Lr_k^{(\ell)})+o(1)\\
  &\ra& 0
  \ena
   as $k\ra\infty$ first then $L\ra\infty$,  that
  \bna
  \int_{\mathbb{B}_{r_k/L}(\widetilde{x}_{i,k})}\varphi_{k}(r_k)\widetilde{f}_k(y,\widetilde{u}_k(y))dy
  \leq N_k(s_k^{(\ell)},r_k/L)+o(1)\ra 0
  \ena
  as $k\ra \infty$ first, then $L\ra\infty$, and that
  $$\int_{\mathbb{B}_{r_k/L}(\widetilde{x}_{i,k})}\varphi_{k}(r_k)\widetilde{u}_k(y)dy\leq
  \int_{\mathbb{B}_{r_k/L}(\widetilde{x}_{i,k})}\widetilde{u}_k^2(y)dy+o(1)\ra 0$$
  as $k\ra\infty$. Therefore we conclude
  \bea\nonumber
  \lim_{L\ra\infty}\limsup_{k\ra\infty}\le|\int_{\mathbb{B}_{1/L}(0)}-\Delta\eta_kdx\ri|&\leq&
  \lim_{L\ra\infty}\limsup_{k\ra\infty}\int_{\mathbb{B}_{r_k/L}(\widetilde{x}_{i,k}))}
  \varphi_{k}(r_k)\widetilde{f}_k(y,\widetilde{u}_k(y))dy\\\nonumber
  &&+\lim_{L\ra\infty}\limsup_{k\ra\infty}\int_{\mathbb{B}_{r_k/L}(\widetilde{x}_{i,k}))}
  \varphi_{k}(r_k)\widetilde{\tau}_k(y)\widetilde{u}_k(y)dy\\\label{(21)}
  &=&0.
  \eea
  Let $\zeta_k$ be a sequence of solution to the equation
  \be\label{(22)}\le\{\begin{array}{lll}
  -\Delta_{\mathbb{R}^2}\zeta_k(x)=e^{\psi_i(\widetilde{x}_{i,k}+r_kx)}\varphi_{k}(r_k)r_k^2
  \widetilde{f}_k(\widetilde{x}_{i,k}+r_kx,v_k(x))\quad{\rm in}\quad
  \mathbb{B}_1(0)\\[1.5ex]\zeta_k=\eta_k\quad{\rm on}\quad
  \p\mathbb{B}_1(0).
  \end{array}\ri.\ee
  Then in view of (\ref{(14)}), $\eta_k-\zeta_k$ satisfies
  \be\label{(23)}\le\{\begin{array}{lll}
  -\Delta_{\mathbb{R}^2}(\eta_k-\zeta_k)(x)=-e^{\psi_i(\widetilde{x}_{i,k}+r_k
  x)}\varphi_{k}(r_k)r_k^2\widetilde{\tau}_k(\widetilde{x}_{i,k}+r_kx)
  v_k(x)\quad{\rm in}\quad
  \mathbb{B}_1(0)\\[1.5ex]\eta_k-\zeta_k=0\quad{\rm on}\quad
  \p\mathbb{B}_1(0).
  \end{array}\ri.\ee
  Since $u_k$ is bounded in $L^p(\Sigma)$ for any $p>1$, applying elliptic estimates to
  (\ref{(23)}), we get
  $$\|\eta_k-\zeta_k\|_{L^\infty(\mathbb{B}_1(0))}\leq C$$
  for some constant $C$. By (\ref{(18)}), $\eta_k$ is uniformly
  bounded on $\p\mathbb{B}_1(0)$. In view of (\ref{(22)}), the
  maximum principle implies that there exists some constant $C$
  such that
  $$\zeta_k(x)\geq -C\quad{\rm for\,\,all}\quad x\in \mathbb{B}_1(0).$$
  Hence
  \be\label{(24)}\eta_k(x)\geq -C\quad{\rm for\,\,all}\quad x\in \mathbb{B}_1(0).\ee
  By (\ref{cla}), $\varphi_{k}(r_k)\leq v_k(x)+C$ for all $x\in\mathbb{B}_{1/L}(0)$ and $L>1$.  Note that
  \bea\nonumber
  \varphi_{k}(r_k)r_k^2\widetilde{f}_k(\widetilde{x}_{i,k}+r_kx,v_k(x))&=&
  \varphi_{k}(r_k)r_k^2\theta_{k}(r_k)\f{\widetilde{f}_k(\widetilde{x}_{i,k}+r_kx,v_k(x))}
  {\theta_{k}(r_k)}\\\label{60'}
  &=&(\alpha_0+o(1))e^{(1+o(1))(v_k^2(x)-\varphi_{k}^2(r_k))}.
  \eea
  Using the inequality $a^2-b^2\geq 2b(a-b)$, $a,b\geq 0$, we get
  $v_k^2(x)-\varphi_{i,k}^2(r_k)\geq 2\eta_k(x)$ for all $x\in\mathbb{B}_1(0)$.
  Then (\ref{60'}) leads to
  \be\label{(25)}\int_{\mathbb{B}_{1/L}(0)}e^{\eta_k}dx\leq \f{2}{\alpha_0}\int_{\mathbb{B}_{1/L}(0)}
  \varphi_{k}(r_k)r_k^2\widetilde{f}_k(\widetilde{x}_{i,k}+r_kx,v_k(x))dx\ee
  for sufficiently large $k$.
  Combining (\ref{(14)}), (\ref{(21)}), (\ref{(24)}) and (\ref{(25)}), we obtain
  \be\label{(26)}\lim_{L\ra\infty}\lim_{k\ra\infty}\int_{\mathbb{B}_{1/L}(0)}\eta_kdx=0.\ee
  For any $\varphi\in C_0^\infty(\mathbb{R}^2)$, integration by
  parts gives
  \bea\nonumber
  \int_{\mathbb{R}^2}\eta\Delta\varphi dx&=&\lim_{L\ra\infty}\int_{\mathbb{R}^2\setminus
  \mathbb{B}_{1/L}(0)}\eta\Delta\varphi dx\\\label{(27)}
  &=&\lim_{L\ra\infty}\le(-\int_{\p\mathbb{B}_{1/L}(0)}\eta\p_\nu\varphi d\sigma
  +\int_{\p\mathbb{B}_{1/L}(0)}\varphi\p_\nu\eta d\sigma+
  \int_{\mathbb{R}^2\setminus
  \mathbb{B}_{1/L}(0)}\varphi\Delta\eta dx\ri).
  \eea
  It is clear that
  \bea\nonumber
  \int_{\p\mathbb{B}_{1/L}(0)}\eta\p_\nu\varphi d\sigma&=&\lim_{k\ra\infty}
  \int_{\p\mathbb{B}_{1/L}(0)}\eta_k\p_\nu\varphi d\sigma\\\nonumber
  &=&\lim_{k\ra\infty}\le(\int_{\mathbb{B}_{1/L}(0)}\eta_k\Delta\varphi dx+
   \int_{\mathbb{B}_{1/L}(0)}\nabla\eta_k\nabla\varphi
   dx\ri)\\\label{(28)}
   &=&\lim_{k\ra\infty}\le(\int_{\mathbb{B}_{1/L}(0)}\eta_k\Delta\varphi dx+
   \int_{\p\mathbb{B}_{1/L}(0)}\varphi\p_\nu\eta_k d\sigma-
   \int_{\mathbb{B}_{1/L}(0)}\varphi\Delta\eta_k dx\ri).
  \eea
  Moreover, by Proposition 4.1 and (\ref{cla}), there exists some
  constant $C$ such that
  $$|\nabla\eta_k(x)|=\varphi_k(r_k)r_k|\nabla u_k(\widetilde{x}_{i,k}+r_k x)|\leq {C}/{|x|}$$
  for all $x\in\mathbb{B}_{1/L}(0)$. This together with (\ref{(21)})
  leads to
  \bna
  \lim_{L\ra\infty}\lim_{k\ra\infty}\int_{\p\mathbb{B}_{1/L}(0)}\varphi\p_\nu\eta_k
  d\sigma&=&\varphi(0)\lim_{L\ra\infty}\lim_{k\ra\infty}\int_{\p\mathbb{B}_{1/L}(0)}\p_\nu\eta_k
  d\sigma\\&=&\varphi(0)\lim_{L\ra\infty}\lim_{k\ra\infty}\int_{\mathbb{B}_{1/L}(0)}\Delta\eta_kdx\\
  &=&0.
  \ena
  As a consequence
  \be\label{(29)}\lim_{L\ra\infty}\int_{\p\mathbb{B}_{1/L}(0)}\varphi\p_\nu\eta
  d\sigma=\lim_{L\ra\infty}\lim_{k\ra\infty}\int_{\p\mathbb{B}_{1/L}(0)}\varphi\p_\nu\eta_k
  d\sigma=0.\ee
  Inserting (\ref{(19)}), (\ref{(21)}), (\ref{(26)}), (\ref{(28)}) and (\ref{(29)}) into
  (\ref{(27)}), we obtain
  $$-\int_{\mathbb{R}^2}\eta\Delta\varphi dx=\lim_{L\ra\infty}\int_{\mathbb{R}^2\setminus\mathbb{B}_{1/L}(0)}
 \alpha_0e^{2\eta}\varphi dx= \int_{\mathbb{R}^2}\alpha_0e^{2\eta}\varphi dx.$$
 Therefore $\eta$ is a distributional solution to the equation
 $$-\Delta_{\mathbb{R}^2}\eta=\alpha_0 e^{2\eta}\quad{\rm in}\quad \mathbb{R}^2.$$
 By the regularity theory for elliptic equations, see for example (\cite{Aubin}, Chapter 2),
 $\eta\in
 C^\infty(\mathbb{R}^2)$. By a result of Chen-Li \cite{C-L},
 $$\eta(x)=\log\f{2}{1+|x|^2}-\log\sqrt{\alpha_0},$$
 and thus
 $$\int_{\mathbb{R}^2}e^{2\eta}dx=\f
{4\pi}{\alpha_0}.$$ This completes the proof of the lemma.
$\hfill\Box$\\

It follows from Lemma 5.6 that
$$
 \lim_{L\ra\infty}\lim_{k\ra\infty}N_k(r_k^{(\ell+1)}/L,Lr_k^{(\ell+1)})
 =\alpha_0\int_{\mathbb{R}^2}e^{2\eta^{(\ell+1)}}dx=4\pi.
$$
 This together with (\ref{(8)}) gives
 $$\lim_{L\ra\infty}\lim_{k\ra\infty}N_k(s_k^{(\ell)},Lr_k^{(\ell+1)})=4\pi.$$
 By the inductive hypothesis $(\mathcal{H}_{\ell,3})$,
 \bna
 \lim_{L\ra\infty}\lim_{k\ra\infty}\Lambda_k(Lr_k^{(\ell+1)})&=&\lim_{L\ra\infty}\lim_{k\ra\infty}
 \le(\Lambda_k(s_k^{(\ell)})+N_k(s_k^{(\ell)},Lr_k^{(\ell+1)})\ri)\\
 &=&4\pi(\ell+1).
 \ena

 Now we set
 $w_k^{(\ell+1)}(x)=\varphi_k(r_k^{(\ell+1)})(\widetilde{u}_k(x)-\varphi_k(r_k^{(\ell+1)}))$.
 Similar to Lemma 5.2, we have\\

 \noindent{\bf Lemma 5.7} {\it For any $\epsilon>0$, let $T_k^{(\ell+1)}=T_k^{(\ell+1)}(\epsilon)>r_k^{(\ell+1)}$ be the minimal
 number such that
 $\varphi_{k}(T_k^{(\ell+1)})=\epsilon\varphi_{k}(r_k^{(\ell+1)})$.
 Then  $r_k^{(\ell+1)}/T_k^{(\ell+1)}\ra 0$ as $k\ra\infty$.
 Moreover, for any $b<2$ and sufficiently large $k$, $L$, there holds
 $$\overline{w_k^{(\ell+1)}}(r)\leq b\log\f{r_k^{(\ell+1)}}{r}+C\quad {\rm for\,\,all}\quad Lr_k^{(\ell+1)}\leq
 r\leq T_k^{(\ell+1)},$$
 where $C$ is a constant depending only on $\alpha_0$ and
 $(\Sigma,g)$, and we have
 $$\lim_{k\ra\infty}N_k(s_k^{(\ell)},T_k^{(\ell+1)})=4\pi.$$}
 {\it Proof.} Since the proof is completely analogous to that of Lemma
 4.2, except that instead of Proposition 3.1 we shall use Lemma
 5.6, the details are omitted here. $\hfill\Box$\\

 For suitable $s_k^{(\ell+1)}=T_k^{(\ell+1)}(\epsilon_k)$, where $\epsilon_k\searrow
 0$ is chosen such that $u_k(s_k^{(\ell+1)})\ra\infty$ as
 $k\ra\infty$ and $r_k^{(\ell+1)}/s_k^{(\ell+1)}\ra 0$ as
 $k\ra\infty$. Moreover
 $$\lim_{k\ra\infty}\Lambda_k(s_k^{(\ell+1)})=4\pi(\ell+1),$$
 and
 $$\lim_{L\ra\infty}\lim_{k\ra\infty}N_k(Lr_k^{(\ell+1)},s_k^{(\ell+1)})=0.$$
 By the definition of $s_k^{(\ell+1)}$,
 $$\lim_{k\ra\infty}\f{\varphi_{k}(s_k^{(\ell+1)})}{\varphi_{k}(Lr_k^{(\ell+1)})}=0
 \quad{\rm for\,\,any}\quad L>0.$$
 Hence
 $(\mathcal{H}_{\ell+1})$ holds. This completes Step 2, and thus the proof of Proposition 5.1.

 \subsection{Quantization for non-simple blow-up points\\}

 In this subsection, we shall prove a quantization result for non-simple blow-up points.
 We assume that $x_i^\ast$ is a non-simple blow-up point of order $m$,
 namely there exists a subset $\{i_1,\cdots,i_m\}\subset\{1,\cdots,N\}$ such that $d_g(x_i^\ast,x_\ell^\ast)=0$
 for all $\ell\in\{i_1,\cdots,i_m\}$ and $d_g(x_j^\ast,x_i^\ast)>0$
 for all $j\in \{1,\cdots,N\}\setminus\{i_1,\cdots,i_m\}$.
 In particular, $i\in\{i_1,\cdots,i_m\}$.
  Take an isothermal coordinate system $(U,\phi;\{x^1,x^2\})$ near $x_i^\ast$,
  where $U\subset\Sigma$ is a neighborhood  of $x_i^\ast$ such that $x_j^\ast\not\in \overline{U}$,
  the closure of $U$ for all $j\in\{1,\cdots,N\}\setminus\{i_1,\cdots,i_m\}$, $\phi:U\ra\Omega\subset\mathbb{R}^2$ is a diffeomorphism
  with $\phi(x_i^\ast)=(0,0)$. We can find some $\delta>0$ such that $\mathbb{B}_{2\delta}(0)\subset
  \Omega$. In this coordinate system, the metric $g=e^{\psi}(d{x^1}^2+d{x_2}^2)$ for some smooth function $\psi:\Omega\ra\mathbb{R}$
  with $\psi(0,0)=0$. We shall prove the following result. \\

  \noindent{\bf Proposition 5.8} {\it Let $u_k$, $u_\infty$, $\tau_k$, $\tau_\infty$, $x_{i,k}$ and $x_i^\ast$ be as in Proposition 3.1.
  Suppose that $x_i^\ast$ is a
  non-simple blow-up point of order $m$ as above. Then up to a subsequence, there exists some positive integer $I$
  such that
  \be\label{identity2}
  \lim_{k\ra\infty}\int_{U}
  (|\nabla_gu_k|^2+\tau_ku_k^2)dv_g=
  \int_{U}(|\nabla_gu_\infty|^2+\tau_\infty
  u_\infty^2)dv_g+4\pi I,
  \ee
  where $U$ is a neighborhood of $x_i^\ast$ chosen as above.}\\

   Similarly as
  before we denote $\widetilde{x}_{j,k}=\phi(x_{j,k})$ for
  $j\in\{i_1,\cdots,i_m\}$, $\widetilde{u}_k=u_k\circ\phi^{-1}$,
  $\widetilde{\tau}_k=\tau_k\circ\phi^{-1}$, and
  $\widetilde{f}_k(x,\widetilde{u}_k(x))=f(\phi^{-1}(x),u_k(\phi^{-1}(x)))$.
  Let $\varphi_k=\varphi_k^{(i)}$, $\Lambda_k=\Lambda_k^{(i)}$ and $N_k=N_k^{(i)}$ be as defined in
  (\ref{mean-val}), (\ref{energy}) and (\ref{neck}) respectively.
  The proof of Proposition 5.8 will be  divided into several steps below.\\

  {\it Step 1. Blow-up analysis at the scale $o(\rho_k)$, where
  $$\rho_k=\rho_k^{(i)}=\f{1}{2}\inf_{j\in\{i_1,\cdots,i_m\}\setminus\{i\}}|\widetilde{x}_{j,k}-\widetilde{x}_{i,k}|.$$ }

  By Proposition 3.1 we have
  $\lim_{L\ra\infty}\lim_{k\ra\infty}\Lambda_k(Lr_{i,k})=4\pi$. Let
  $r_k^{(1)}=r_{i,k}$.
   We distinguish the following two
 cases to proceed.\\

 \noindent{\bf Case 1} {\it there exists some $0<\epsilon_0<1$ such that for all $t\in[r_k^{(1)},\rho_k]$ there holds
 $\varphi_k(t)\geq \epsilon_0\varphi_k(r_k^{(1)})$};\\[1.2ex]
 {\bf Case 2} {\it for any $\epsilon>0$ there exists a minimal $T_k=T_k(\epsilon)\in[r_k^{(1)},\rho_k]$
 such that $\varphi_k(T_k)=\epsilon\varphi_k(r_k^{(1)})$.}\\

 In Case 1, the decay estimate that we established in Lemma 5.2
 remains valid on $[r_k^{(1)},\rho_k]$. Moreover
 $$\lim_{L\ra\infty}\lim_{k\ra\infty}\Lambda_k(s_k)=4\pi$$
 for any sequence $s_k$ satisfying $s_k/\rho_k\ra 0$ and $s_k/r_k^{(1)}\ra\infty$
 as $k\ra\infty$. The concentration analysis at scales up to
 $o(\rho_k)$ is complete.

 In Case 2, as before we can find numbers
 $s_k^{(1)}<\rho_k$ with $\varphi_k(s_k^{(1)})\ra\infty$ as
 $k\ra\infty$,
 $\Lambda_k(s_k^{(1)})\ra 4\pi$ as $k\ra\infty$,
 and
 $\varphi_k(s_k^{(1)})/\varphi_k(Lr_k^{(1)})\ra 0$ for any $L\geq 1$ as $k\ra\infty$.
 We proceed by iteration up to some maximal index $\ell_0\geq 1$
 where either Case 1 or (\ref{caseA}) holds with final radii
 $r_k^{(\ell_0)},s_k^{(\ell_0)}$, respectively. Hence
 \be\label{(72)}\lim_{k\ra\infty}\Lambda_k(s_k^{(\ell_0)})=4\pi\ell_0,
 \quad \lim_{k\ra\infty}{\varphi_k(s_k^{(\ell_0)})}/{\varphi_k(Lr_k^{(\ell_0)})}=0,\,\,\forall L\geq 1\ee
 and
 \be\label{73}\lim_{k\ra\infty}N_k(s_k^{(\ell_0)},t_k)=0\,\,\,{\rm
 for\,\,any\,\,sequence}\,\,\,
 t_k=o(\rho_k).\ee
 This leads to
 \be\label{74}\lim_{L\ra\infty}\lim_{k\ra\infty}N_k(s_k^{(\ell_0)},\rho_k/L)=0.\ee
 For otherwise, we can find some $\mu_0>0$ such that up to a
 subsequence
 $$\lim_{k\ra\infty}N_k(s_k^{(\ell_0)},\rho_k)\geq\lim_{k\ra\infty}N_k(s_k^{(\ell_0)},\rho_k/L)\geq\mu_0$$
 for all $L\geq 1$. Take $t_k^\prime\in(s_k^{(\ell_0)},\rho_k)$ such that
 \be\label{75}0<\lim_{k\ra\infty}N_k(s_k^{(\ell_0)},t_k^\prime)<\f{1}{2C_0},\ee
 where $C_0$ is a constant as in Lemma 5.3. Then by Lemma 5.4 we
 have
 $$\lim_{L\ra\infty}\lim_{k\ra\infty}N_k(s_k^{(\ell_0)},t_k^\prime/L)=0.$$
 In view of (\ref{73}) and (\ref{75}), there exists some $\nu_0>0$
 such that up to a subsequence, $t_k^\prime\geq \nu_0\rho_k$ for all $k$. This immediately implies
 (\ref{74}) and completes Step
 1.\\

To proceed, we introduce several terminologies concerning the
classification of blow-up points near $x_i^\ast$.
Define a set
$$\mathcal{X}=\mathcal{X}^{(i)}=\{x_{i_1,k},\cdots,x_{i_m,k}\},$$
where each $x_{j,k}$, $j\in\{i_1,\cdots,i_m\}$, denotes a sequence $(x_{j,k})$. In the sequel we do not distinguish sequences $(x_{j,k})$ and
points $x_{j,k}$.
Let $t_k>0$ be a bounded sequence. For any $j\in\{i_1,\cdots,i_m\}$, we define a $t_k$-equivalent class associated to the
sequence $x_{j,k}$ by
$$[x_{j,k}]_{t_k}:=\le\{x_{\ell,k}:
d_g(x_{\ell,k},x_{j,k})=o(t_k),\,\,\ell\in \{i_1,\cdots,i_m\}\ri\}.$$
The total number of sequences in $[x_{j,k}]_{t_k}$ is called the order of $[x_{j,k}]_{t_k}$.
In particular, the order of $[x_{j,k}]_{\rho_k^{(j)}}$ is exactly one, while the order of $[x_{j,k}]_{\delta}$ is $m$. Actually
we have $[x_{j,k}]_{\delta}=\mathcal{X}$.
 Moreover,
if $x_{\ell,k}\in [x_{j,k}]_{t_k}$, then
$x_{j,k}\in[x_{\ell,k}]_{t_k}$. Also, if
$[x_{j,k}]_{t_k}\cap[x_{\ell,k}]_{t_k}\not=\varnothing$, then
$[x_{j,k}]_{t_k}=[x_{\ell,k}]_{t_k}$. Hence every subset of $\mathcal{X}$
can be divided into several $t_k$-equivalent classes, any two of
which have no intersection.

For any $1\leq \ell<m$, we say that the property
$(\mathcal{A}_\ell)$ holds for some $t_k$-equivalent class
$[x_{j,k}]_{t_k}$ of order $\ell$, if
  either $(a)$
  there exist $r_k>0$ and integer $I^{(j)}$ such that for some
  $\epsilon_0>0$ and
 all $t\in[r_k,t_k]$ there holds $\varphi_k^{(j)}(t)\geq
 \epsilon_0\varphi_k^{(j)}(r_k)$, $\Lambda_k^{(j)}(Lr_k)\ra 4\pi
 I^{(j)}$ and $N_k^{(j)}(Lr_k,t_k/L)\ra 0$ as $k\ra\infty$ first, and then $L\ra\infty$;
   or
 $(b)$ there exist sequences $r_k<s_k<t_k$ and an integer $I^{(j)}$ such that
 $\varphi_k^{(j)}(s_k)/\varphi_k^{(j)}(Lr_k)\ra 0$  as $k\ra\infty$ for any $L\geq
 1$, $\Lambda_k^{(j)}(t_k/L)\ra 4\pi I^{(j)}$ and $N_k^{(j)}(s_k,t_k/L)\ra 0$ as
 $k\ra\infty$ first, and then $L\ra\infty$.
 While we say that the property $(\mathcal{A}_m)$ holds, if there exits some $j\in\{i_1,\cdots,i_m\}$ and
 integer $I^{(j)}$ such that $\Lambda_k^{(j)}(\delta/L)\ra 4\pi I^{(j)}$ as $k\ra\infty$ first, and then $L\ra\infty$.

 According to Proposition 5.1, when $m=1$, $(\mathcal{A}_1)$
 holds.
 When $m>1$, we let $\rho_{k,0}=\rho_k$ and $\rho_{k,j}$ $(1\leq j\leq m-1)$ be defined as in (\ref{5.97}) and (\ref{5.100}) below.
 It follows from Step 1 that $(\mathcal{A}_1)$ holds for any $t_k$-equivalent class of order one, where
 \be\label{tk}t_k\in\{\rho_{k,0},\cdots,\rho_{k,m-1}\}.\ee
 We now we make an induction procedure on both orders of $t_k$-equivalent class and $m$. Suppose that for some integer $\nu\geq 1$, when
 $m=\nu$, the property $(\mathcal{A}_\nu)$ holds; while when $m>\nu$, the property $(\mathcal{A}_\ell)$
 holds for any  $t_k$-equivalent class of order $1\leq\ell\leq \nu$, where $t_k$ is as in (\ref{tk}). We shall prove the following:
 When $m=\nu+1$, the property $(\mathcal{A}_{\nu+1})$ holds; When $m>\nu+1$, the property $(\mathcal{A}_\ell)$ holds for any
 $t_k$-equivalent class of order $1\leq\ell\leq \nu+1$, where $t_k$ is as in (\ref{tk}).
 Assuming this induction argument is complete, we conclude that $(\mathcal{A}_m)$ holds for any integer $m$.
 It is easy to see that (\ref{identity2})
 follows immediately from  $(\mathcal{A}_m)$ and the fact that
 $u_k\ra u_\infty$ in $C^1_{\rm loc}(\Sigma\setminus\{x_1^\ast,\cdots,x_i^\ast\})$.\\

 In the next two steps, we shall prove that $(\mathcal{A}_m)$ holds for $m=\nu+1$. In Step 4, we shall prove that
 $(\mathcal{A}_\ell)$ holds for any
 $t_k$-equivalent class of order $1\leq\ell\leq \nu+1$, where $t_k$ is as in (\ref{tk}).
 \\

 {\it Step 2. Blow-up analysis at the scale $\rho_k$.}\\

 Let $m=\nu+1$. Now we turn to carry out blow-up analysis at the scale $\rho_k$ near
 $\widetilde{x}_{i,k}$. We first
 assume that for some $L\geq 1$ there exists some sequence $(x_k)$
 such that $\rho_k/L\leq R_k(x_k)\leq|x_k-\widetilde{x}_{i,k}|\leq L\rho_k$ and
 \be\label{76}|x_k-\widetilde{x}_{i,k}|^2\widetilde{u}_k(x_k)\widetilde{f}_k(x_k,\widetilde{u}_k(x_k))\geq \nu_0>0.\ee
 By Proposition 4.1 we may assume that $|x_k-\widetilde{x}_{i,k}|=\rho_k$. The
 following estimate is important for our subsequent analysis.\\

 \noindent{\bf Lemma 5.9} {\it Assuming (\ref{76}), we have $\varphi_k(\rho_k)/\varphi_k(r_k^{(\ell_0)})\ra
 0$
 as $k\ra \infty$.}\\

 \noindent{\it Proof.} If we suppose that there exists some
 $\epsilon_0>0$ such that $\varphi_k(\rho_k)\geq
 \epsilon_0\varphi_k(r_k^{(\ell_0)})$, then we set
 $$w_k(x)=\varphi_k(r_k^{(\ell_0)})(\widetilde{u}_k(x)-\varphi_k(r_k^{(\ell_0)})),\quad x\in\Omega.$$
 Similar to Lemma 5.2, there holds for any $b<2$
 \be\label{wdecay}\overline{w}_k(r)\leq
 b\log\f{r_k^{(\ell_0)}}{r}+C\ee
 for all $r\in [r_k^{(\ell_0)},\rho_k]$. Let $\theta_k$ be as defined in (\ref{theta}). By (H5) and $(iii)$
 of Proposition 3.1, we find some uniform constant
 $C$ such that
 \be\label{rphi}r_k^{(\ell_0)}\varphi_k(r_k^{(\ell_0)})\theta_k(r_k^{(\ell_0)})\leq C.\ee
   Hence we obtain
 \bea\nonumber
 |x_k-\widetilde{x}_{i,k}|^2u_k(x_k)\widetilde{f}_k(x_k,\widetilde{u}_k(x_k))&\leq&
 C\rho_k^2\varphi_k(\rho_k)\theta_k(\rho_k)\\
 \nonumber&=&C(r_k^{(\ell_0)})^2\varphi_k(r_k^{(\ell_0)})\theta_k(r_k^{(\ell_0)})\le(\f{\rho_k}{r_k^{(\ell_0)}}\ri)^2
 \f{\varphi_k(\rho_k)}{\varphi_k(r_k^{(\ell_0)})}
 \f{\theta_k(\rho_k)}{\theta_k(r_k^{(\ell_0)})}\\
 \nonumber&\leq&
 C\le({\rho_k}/{r_k^{(\ell_0)}}\ri)^2e^{(1+o(1))(\varphi_k^2(\rho_k)-\varphi_k^2(r_k^{(\ell_0)}))}\\
 \nonumber&\leq&C\le({\rho_k}/{r_k^{(\ell_0)}}\ri)^2e^{(1+o(1))(1+\epsilon_0)\overline{w}_k(\rho_k)}\\
 &\leq&C\le(\rho_k/r_k^{(\ell_0)}\ri)^{2-(1+o(1))(1+\epsilon_0)b}\ra
 0 \label{c-t}
 \eea
 as $k\ra\infty$, if we choose $b<2$ such that $(1+\epsilon_0)b>2$.
 Here the first inequality follows from Proposition 4.1, the second
 one follows from (H4), (H5) and (\ref{rphi}), while the
 third one is a consequence of our assumption $\varphi_k(\rho_k)\geq
 \epsilon_0\varphi_k(r_k^{(\ell_0)})$, and the last one is implied by (\ref{wdecay}).
 The contradiction between (\ref{c-t})
 and (\ref{76}) ends the proof of the lemma. $\hfill\Box$\\

 Lemma 5.9 implies that for any $\epsilon>0$ there exists $T_k\in
 [r_k^{(\ell_0)},\rho_k]$ such that
 $\varphi_k(T_k)=\epsilon\varphi_k(r_k^{(\ell_0)})$. Hence at scales
 up to order $o(\rho_k)$ we end up with (\ref{caseA}), where $\ell$ is replaced by $\ell_0$. The desired quantization
 result at the scale $\rho_k$ then is a consequence of the following
 result.\\

 \noindent{\bf Lemma 5.10} {\it Assuming (\ref{76}), then up to a
 subsequence we can find some $\alpha_0\geq \nu_0$ such that
 \be\label{761}\lim_{k\ra\infty}|x_k-\widetilde{x}_{i,k}|^2\widetilde{u}_k(x_k)\widetilde{f}_k(x_k,\widetilde{u}_k(x_k))=\alpha_0.\ee
 Moreover there exist a finite set $\mathcal{S}_\infty\subset\mathbb{R}^2$ such that
 $$\eta_k(x)=\widetilde{u}_k(x_k)(\widetilde{u}_k(\widetilde{x}_{i,k}+\rho_kx)-\widetilde{u}_k(x_k))\ra\eta(x)=\log\f{2}{\sqrt{\alpha_0}(1+|x|^2)}$$
 in $C^1_{\rm loc}(\mathbb{R}^2\setminus\mathcal{S}_\infty)$ as
 $k\ra\infty$.}\\

 \noindent{\it Proof.} It is obvious that (\ref{761}) holds for some $\alpha_0\geq \nu_0>0$.  Define
 $$v_k(y)=\widetilde{u}_k(\widetilde{x}_{i,k}+\rho_ky)$$
 for $y\in\Omega_k=\{y\in\mathbb{R}^2:\widetilde{x}_{i,k}+\rho_ky\in
 \Omega\}$. Let
 $$y_{j,k}=\f{\widetilde{x}_{j,k}-\widetilde{x}_{i,k}}{\rho_k}$$
 and
 $$\mathcal{S}_k=\mathcal{S}^{(i)}_{k}=\le\{y_{j,k}: j=i_1,\cdots,i_{\nu+1}\ri\}.$$
 Without loss of generality we assume either $|y_{j,k}|\ra\infty$ or
 $y_{j,k}\ra y_j$, $j=i_1,\cdots,i_{\nu+1}$, and we let
 $\mathcal{S}_\infty=\mathcal{S}_{\infty}^{(i)}$ be the set of accumulation
 points of $\mathcal{S}_k$. Also we let
 $$y_{0,k}=\f{x_k-\widetilde{x}_{i,k}}{\rho_k}$$
 be the scaled points of $x_k$ for which (\ref{76}) holds and which
 satisfy $|y_{0,k}|=1$. Moreover we can assume $y_{0,k}\ra y_0$ as
 $k\ra\infty$.

 Since $\widetilde{u}_k(x_k)\ra\infty$ by (\ref{76}) and
 $\mathcal{S}_\infty$ is a finite set, we have by using Proposition 4.1 and a standard covering argument that
 \be\label{77}v_k-\widetilde{u}_k(x_k)\ra 0\,\,{\rm locally\,\, uniformly\,\,on}\,\,\mathbb{R}^2\setminus\mathcal{S}_\infty\ee
 as $k\ra\infty$. Using the same argument as in the proof of
 Lemma 5.6, we obtain
 $$\eta_k\ra\eta\quad{\rm in}\quad C^1_{\rm
 loc}(\mathbb{R}^2\setminus\mathcal{S}_\infty),$$
 where $\eta\in C^\infty(\mathbb{R}^2\setminus\mathcal{S}_\infty)$
 satisfies the equation
 $$-\Delta_{\mathbb{R}^2}\eta=\alpha_0 e^{2\eta}\quad{\rm in}\quad
 \mathbb{R}^2\setminus\mathcal{S}_\infty.$$
 It follows from (\ref{77}) that $v_k/\widetilde{u}_k(x_k)\ra 1$ locally uniformly on $\mathbb{R}^2\setminus
 \mathcal{S}_\infty$. For any $L\geq 1$ we write
 $$K_L=\mathbb{B}_L(0)\setminus(\cup_{y_j\in\mathcal{S}_\infty}\mathbb{B}_{\delta/L}(y_j)).$$
 Combining (H4), (H5), (\ref{bdd}) and (\ref{77}), we can estimate
 \bna
 \int_{\mathbb{R}^2}e^{2\eta}dx&\leq&
 \lim_{L\ra\infty}\lim_{k\ra\infty}\int_{K_L}\f{v_k(x)}{\widetilde{u}_k(x_k)}e^{(1+o(1))\eta_k(1+\f{v_k(x)}{\widetilde{u}_k(x_k)})}dx\\
 &=&\lim_{L\ra\infty}\lim_{k\ra\infty}\int_{K_L}\f{\widetilde{u}_k(\widetilde{x}_{i,k}+\rho_kx)
 \widetilde{f}_k(\widetilde{x}_{i,k}+\rho_kx,\widetilde{u}_k(\widetilde{x}_{i,k}+\rho_kx))}
 {\widetilde{u}_k(x_k)\widetilde{f}_k(x_k,\widetilde{u}_k(x_k))}dx\\
 &\leq&\f{C}{\nu_0}\limsup_{k\ra\infty}\int_\Sigma
 u_kf_k(x,u_k)dv_g\leq \f{C}{\nu_0}.
 \ena

 Since $y_{j,k}\ra y_j$ as $k\ra\infty$, we can take sufficiently large $L$ and $k$ such that
  $\mathbb{B}_{1/L}(y_j)\subset\mathbb{B}_{2/L}(y_{j,k})$ and
  $\mathbb{B}_{2/L}(y_{j,k})\cap \mathbb{B}_{2/L}(y_{\alpha,k})=\varnothing$ for any
  $\alpha\not=j$. Moreover let $\ell$ be the order of the $\rho_k$-equivalent class $[x_{j,k}]_{\rho_k}$. Clearly
  $\ell\leq \nu$.
  By our inductive assumption,
  $(\mathcal{A}_\ell)$ holds for $[x_{j,k}]_{\rho_k}$. Noting that Lemma 5.9 excludes the possibility of Case 1
  with $r_k^{(1)}$ replaced by $r_k^{(\ell_0)}$, we can find
  sequences $r_k^{(j)}<s_k^{(j)}$ such that
  \be\label{bz}\lim_{k\ra\infty}\varphi_k(s_k^{(j)})/\varphi_k(Lr_k^{(j)})=0,\quad\forall L\geq 1.\ee
  and
  \be\label{nz}\lim_{L\ra\infty}\lim_{k\ra\infty}N_k^{(j)}(s_k^{(j)},\rho_k/L)=0,\ee
   Note again that $y_{j,k}\ra y_j$ as $k\ra\infty$. There exists some constant $C$,
    which may depends on $|y_j|$ but
   not on $k$, such that $|\widetilde{x}_{j,k}-\widetilde{x}_{i,k}|\leq C\rho_k$.
   For any
   ${x}_{\alpha,k}\not\in[{x}_{j,k}]_{\rho_k}$,
   we can take some large $L_0$ such that
   $|\widetilde{x}_{j,k}-\widetilde{x}_{\alpha,k}|\geq \rho_k/(2L_0)$ for all sufficiently large $k$. Recalling that
   $|x_k-\widetilde{x}_{i,k}|=\rho_k$ and applying Proposition
   4.1, we obtain
   $$\widetilde{u}_k(x_k)\leq \inf_{\p\mathbb{B}_{2\rho_k/L_0}(\widetilde{x}_{j,k})}\widetilde{u}_k+C$$
   for some uniform constant $C$. As in the proof of Lemma 5.2, we can
   find another uniform constant $C$ such that for all $x\in \mathbb{B}_{2\rho_k/L_0}(\widetilde{x}_{j,k})$
   $$\widetilde{u}_k(x)\geq \inf_{\p\mathbb{B}_{2\rho_k/L_0}(\widetilde{x}_{j,k})}\widetilde{u}_k-C.$$
   These two estimates immediately imply the existence of some uniform constant $C$ such that
   \be\label{ukleq}\widetilde{u}_k(x_k)\leq \widetilde{u}_k(x)+C\,\,{\rm for\,\,all}\,\,x\in\mathbb{B}_{2\rho_k/L}
   (\widetilde{x}_{j,k}),
   \ee provided that $L\geq L_0$. Note that $g=e^\psi(d{x^1}^2+d{x^2}^2)$ for some smooth function $\psi$ with
    $\psi(0,0)=0$. By the equation (\ref{e1}), we have for large $L$
   \bna
   \int_{\mathbb{B}_{1/L}(y_j)}|\Delta_{\mathbb{R}^2}\eta_k|dx&\leq&
   \int_{\mathbb{B}_{1/L}(y_j)}\rho_k^2\widetilde{u}_k(x_k)\widetilde{f}_k\le(\widetilde{x}_{i,k}+\rho_kx,v_k(x)\ri)
   e^{\psi(\widetilde{x}_{i,k}+\rho_kx)}dx\\
   &&+\int_{\mathbb{B}_{1/L}(y_j)}\rho_k^2\widetilde{u}_k(x_k)\widetilde{\tau}_k(\widetilde{x}_{i,k}+\rho_kx)
   v_k(x)e^{\psi(\widetilde{x}_{i,k}+\rho_kx)}dx\\
   &\leq&\int_{\mathbb{B}_{2/L}(y_{j,k})}\rho_k^2\widetilde{u}_k(x_k)\widetilde{f}_k
   \le(\widetilde{x}_{i,k}+\rho_kx,v_k(x)\ri)e^{\psi(\widetilde{x}_{i,k}+\rho_kx)}dx\\
   &&+\int_{\mathbb{B}_{2/L}(y_{j,k})}\rho_k^2\widetilde{u}_k(x_k)\widetilde{\tau}_k(\widetilde{x}_{i,k}+\rho_kx)
   v_k(x)e^{\psi(\widetilde{x}_{i,k}+\rho_kx)}dx\\
   &=&\int_{\mathbb{B}_{2\rho_k/L}(\widetilde{x}_{j,k})}\widetilde{u}_k(x_k)\le(\widetilde{f}_k(y,\widetilde{u}_k(y))
   +\widetilde{\tau}_k(y) \widetilde{u}_k(y)\ri)e^{\psi(y)}dy.
   \ena
   With the help of (\ref{bz})-(\ref{ukleq}) and an obvious analogy to (\ref{cla}), we obtain
   $$\lim_{L\ra\infty}\lim_{k\ra\infty}\int_{\mathbb{B}_{1/L}(y_j)}|\Delta_{\mathbb{R}^2}\eta_k|dx=0,$$
   analogous to (\ref{(21)}).
   In the same way of proving (\ref{(26)}) we get
   $$\lim_{L\ra\infty}\lim_{k\ra\infty}\int_{\mathbb{B}_{1/L}(y_j)}\eta_kdx=0.$$

   In view of (\ref{ukleq}), we can find some uniform constant $C$ such that for all $y\in\p\mathbb{B}_{1/L}(y_j)$
   $$\widetilde{u}_k(x_k)/\widetilde{u}_k(\widetilde{x}_{j,k}+\rho_k y)\leq C,$$
   which together with Proposition 4.1 leads to
   $$
   |y-y_{j,k}||\nabla_{\mathbb{R}^2}\eta_k(y)|=|\widetilde{x}_{i,k}+\rho_k
   x-\widetilde{x}_{j,k}|\widetilde{u}_k(x_k)|\nabla_{\mathbb{R}^2} \widetilde{u}_k(\widetilde{x}_{i,k}+\rho_kx)|\leq C.
   $$
   This gives
   $$|\nabla_{\mathbb{R}^2}\eta_k(y)|\leq \f{C}{|y-y_j|}$$
   for all $y\in\p\mathbb{B}_{1/L}(y_j)$, provided that $k$ is
   sufficiently large. Then we obtain an analogy to (\ref{(29)}),
   namely, for any $\varphi\in C_0^\infty(\mathbb{R}^2)$
   $$\lim_{L\ra\infty}\int_{\p\mathbb{B}_{1/L}(y_j)}\varphi\p_\nu\eta
  d\sigma=\lim_{L\ra\infty}\int_{\p\mathbb{B}_{1/L}(y_j)}\eta\p_\nu\varphi
  d\sigma=0.$$
  This excludes $y_j$ as a singular point of $\eta$ as in Lemma 5.6. Since $y_j$ is any point of $\mathcal{S}_\infty$, we
  conclude that $\eta$ is a smooth solution to the equation
  $$-\Delta_{\mathbb{R}^2}\eta=\alpha_0 e^{2\eta}\quad{\rm in}\quad\mathbb{R}^2.$$
  The remaining part of the conclusions of the lemma follows from a
  result of Chen-Li \cite{C-L}. $\hfill\Box$\\

  Define a set
  \be\label{ALk}A_{L,k}=\le\{x\in \Omega: \rho_k/L\leq R_k(x)\leq|x-\widetilde{x}_{i,k}|\leq
  L\rho_k\ri\}.\ee
  It follows from Proposition 4.1 that ${u_k(x)}/{\widetilde{u}_k(x_k)}\ra 1$
  uniformly in $A_{L,k}$ as $k\ra\infty$. Thus by Lemma 5.10, in case of
  (\ref{76}) there holds
  \bea
  \lim_{L\ra\infty}\lim_{k\ra\infty}\int_{A_{L,k}}\widetilde{u}_k(x)\widetilde{f}_k(x,\widetilde{u}_k(x))dx&=&\alpha_0
  \lim_{L\ra\infty}\lim_{k\ra\infty}\int_{A_{L,k}}\f{\widetilde{u}_k(x)\widetilde{f}_k(x,\widetilde{u}_k(x))}
  {\widetilde{u}_k(x_k)\widetilde{f}_k(x_k,\widetilde{u}_k(x_k))}dx{\nonumber}\\
  &=&\alpha_0\int_{\mathbb{R}^2}e^{2\eta(x)}dx{\nonumber}\\\label{80}
  &=&\alpha_0\f{4\pi}{\alpha_0}=4\pi.
  \eea
  Let
  \be\label{Xk1}X_{k,1}=X_{k,1}^{(i)}=\{\widetilde{x}_{j,k}:\exists C>0 \,\,{\rm such\,\, that}\,\,
  |\widetilde{x}_{j,k}-\widetilde{x}_{i,k}|\leq C\rho_k
  \,\,{\rm for\,\,all}\,\, k\}.\ee
  We can divide $X_{k,1}$ into several $\rho_k$-equivalent classes with
  their orders no more than $\nu$.
  Recalling our inductive assumption $(\mathcal{A}_\ell)$ with $1\leq\ell\leq\nu$ and using (\ref{80}),
  we can find some integer $I$ such that
  $$\lim_{L\ra\infty}\lim_{k\ra\infty}\Lambda_k(L\rho_k)=4\pi(1+I).$$

  On the other hand, if (\ref{76}) does not hold, we have
  \be\label{varn}\lim_{L\ra\infty}\lim_{k\ra\infty}\int_{A_{L,k}}\widetilde{u}_k\widetilde{f}_k(x,\widetilde{u}_k)dx=0.\ee
  The energy estimate at the scale $\rho_k$ again is
  finished. \\

  {\it Step 3. Blow-up analysis at scales exceeding
  $\rho_k$.}\\

  Now we deal with blow-up analysis at scales
  exceeding $\rho_k$ near $\widetilde{x}_{i,k}$. Write
 $$X_{k,0}=\{\widetilde{x}_{i_1,k},\cdots,\widetilde{x}_{i_m,k}\}.$$
 Recalling (\ref{Xk1}), we let
  \be\label{5.97}\rho_{k,1}=\rho_{k,1}^{(i)}=\le\{\begin{array}{lll}\inf\limits
  _{\widetilde{x}_{j,k}\in X_{k,0}\setminus X_{k,1}}\f{|\widetilde{x}_{j,k}-\widetilde{x}_{i,k}|}{2} &{\rm
  if}
  &X_{k,0}\setminus X_{k,1}\not=\varnothing\\[1.5ex]
  \delta, &{\rm if}&X_{k,0}\setminus X_{k,1}=\varnothing.\end{array}\ri.\ee
  From this definition it follows that $\rho_{k,1}/\rho_k\ra \infty$ as $k\ra\infty$. Then,
  using the obvious analogy of Lemma 5.4, either we have
  $$\lim_{L\ra\infty}\lim_{k\ra\infty}N_k(L\rho_k,\rho_{k,1}/L)=0,$$
  and we iterate to the next scale; or there exist a sequence $t_k$
  such that $t_k/\rho_k\ra\infty$, $t_k/\rho_{k,1}\ra 0$ as
  $k\ra\infty$ and up to a subsequence such that
  \be\label{85}P_k(t_k)\geq \nu_0>0\,\,{\rm for\,\,all\,\,large}\,\, k.\ee
  The argument then depends on whether (\ref{76}) or (\ref{varn})
  holds. In case of (\ref{76}), as in Lemma 5.9, the bound (\ref{85})
  and Lemma 5.10 imply that $\varphi_k(t_k)/\varphi_k(\rho_k)\ra 0$ as
  $k\ra \infty$. Then we can argue as in (\ref{caseA}) for
  $r\in[L\rho_k,\rho_{k,1}/L]$ for sufficiently large $L$, and we
  can continue as before to resolve concentrations in this range of
  scales.

  In case of (\ref{varn}) we further need to distinguish whether
  (\ref{caseA}) or Case 1 holds at the final stage of our analysis at
  scales
  $o(\rho_k)$. Recalling that in case of (\ref{caseA}) we have (\ref{(72)}) and
  (\ref{74}), in view of (\ref{varn}) for a suitable sequence of
  numbers $s_{k,1}^{(0)}$ such that $s_{k,1}^{(0)}/\rho_k\ra\infty$,
  $t_k/s_{k,1}^{(0)}\ra\infty$ as $k\ra\infty$ we obtain
  \be\label{sum}\lim_{L\ra\infty}\lim_{k\ra\infty}\le(\Lambda_k(s_{k,1}^{(0)})-
  \sum_{\widetilde{x}_{j,k}\in\widetilde{X}_{k,1}}\Lambda_k^{(j)}(Lr_k^{(\ell_0^{(j)})})\ri)=0,\ee
  where $\Lambda_k^{(j)}(r)$ and $r_k^{(\ell_0^{(j)})}$ are computed
  as above with respect to the blow-up point $x_{j,k}$ and
  $\widetilde{X}_{k,1}$ is the modular set containing all $t_k$-equivalent classes
  of ${X}_{k,1}$, whence the distance between any two points of
  $\widetilde{X}_{k,1}$ is greater than $\widetilde{\nu} \rho_k$ for some
  constant $\widetilde{\nu}>0$. In particular, with such a choice of
  $s_{k,1}^{(0)}$ we find the immediate quantization result
  $$\lim_{k\ra\infty}\Lambda_k(s_{k,1}^{(0)})=4\pi I$$
  for some positive integer $I$. Here again we use the inductive assumption that $(\mathcal{A}_l)$ holds for all
  $\rho_k$-equivalent classes of order $\ell$ with $1\leq\ell\leq \nu$.
  While in Case 1 if we assume there is some  $\epsilon_0>0$ such that
  \be\label{102}\varphi_k(s_{k,1}^{(0)})\geq \epsilon_0
  \varphi_k(Lr_k^{(\ell_0^{(j)})})\ee
  for all $r\in[Lr_{k}^{\ell_0^{(j)}},s_{k,1}^{(0)}]$, then
  as  before we have
  $$\lim_{L\ra\infty}\lim_{k\ra\infty}N_k^{(j)}(Lr_k^{(\ell_0^{(j)})},s_{k,1}^{(0)})=0.$$
  This contradicts (\ref{sum}) since $s_{k,1}^{(0)}/\rho_k\ra\infty$ as $k\ra\infty$ and
  the modular set $\widetilde{X}_{k,1}$ has at
  least two elements. This implies that (\ref{102}) does not hold and
  up to a subsequence there holds for any $L\geq 1$
  $$\lim_{k\ra\infty}\f{\varphi_k(s_{k,1}^{(0)})}{\varphi_k(Lr_k^{(\ell_0^{(j)})})}=0$$
  for all $x_{j,k}\in \widetilde{X}_{k,1}$ where Case 1 holds. Then
  we can continue to resolve concentrations on the range
  $[s_{k,1}^{(0)},\rho_{k,1}/L]$ as before.

  We then proceed by iteration. For $\ell\geq 2$ we
  inductively define the sets
  $$X_{k,\ell}=X_{k,\ell}^{(i)}=
  \le\{\widetilde{x}_{j,k}: \exists C>0\,\,{\rm such\,\, that}\,\, |\widetilde{x}_{j,k}-\widetilde{x}_{i,k}|
  \leq C\rho_{k,\ell-1}
  \,\,{\rm for\,\,all}\,\, k\ri\}$$
  and let
  \be\label{5.100}\rho_{k,\ell}=\rho_{k,\ell}^{(i)}=\le\{\begin{array}{lll}\inf\limits
  _{\widetilde{x}_{j,k}\in X_{k,0}\setminus X_{k,\ell}}\f{|\widetilde{x}_{j,k}-\widetilde{x}_{i,k}|}{2} &{\rm
  if}
  &X_{k,0}\setminus X_{k,\ell}\not=\varnothing\\[1.5ex]
  \delta, &{\rm if}&X_{k,0}\setminus X_{k,\ell}=\varnothing.\end{array}\ri.\ee
  Iteratively carrying out the above analysis at all scales
  $\rho_{k,\ell}$, exhausting all blow-up points $x_{j,k}$, up
  to a subsequence we obtain quantization result for
  $X_{k,0}$. Then Step 3 is finished.\\

  It follows from Step 2 and Step 3 that there exists some integer
  $I$ such that
  \be\label{ast}\lim_{L\ra\infty}\lim_{k\ra\infty}\Lambda_k(\delta/L)=4\pi
  I,\ee
  different analogous to Lemma 5.5. Here and in the sequel, $I$ may denote different integer.
  Hence the property $(\mathcal{A}_m)$ holds when $m=\nu+1$.\\

  {\it Step 4. $(\mathcal{A}_\ell)$ holds for $1\leq\ell\leq\nu+1$ when $m>\nu+1$.}\\

  When $m>\nu+1$, by our inductive assumption, $(\mathcal{A}_\ell)$ holds for all $1\leq\ell\leq \nu$,
  it suffices to prove that
  $(\mathcal{A}_{\nu+1})$ holds for any $t_k$-equivalent class
  $[x_{j,k}]_{t_k}$ of order $\nu+1$, where $j\in\{i_1,\cdots,i_m\}$ and $t_k$ is as in (\ref{tk}).
    This is completely analogous to that $(\mathcal{A}_m)$ holds in the
  case of $m=\nu+1$, which we proved above, except that (\ref{ast}) is replaced by
  $$\lim_{L\ra\infty}\lim_{k\ra\infty}\Lambda_k^{(j)}(t_k/L)=4\pi
  I$$
  for some integer $I$. We omit the details here. This ends Step 4. $\hfill\Box$\\

  Proposition 5.8 follows from the property $(\mathcal{A}_m)$ and the last assertion of Proposition 3.1.

 \section{Completion of the proof of Theorem 1.1}
 In this section, we complete the proof of Theorem 1.1.
 Let $x_{i,k}\ra x_i^\ast$ as $k\ra\infty$, $1\leq i\leq N$, be as in Proposition 3.1.
 In view of possible
 non-simple blow-up points, without loss of generality, we
 may assume for some $q\leq N$, $x_1^\ast,\cdots,x_q^\ast$ are different from each other and
 $x_\ell^\ast\in\{x_1^\ast,\cdots,x_q^\ast\}$ for any $q+1\leq \ell\leq
 N$.
 For any $1\leq i\leq q$, we take an isothermal coordinate system $(U_i,\phi_i;\{x^1,x^2\})$
 near $x_i^\ast$ such that $\phi_i(x_i^\ast)=(0,0)$ and $U_i=\phi_i^{-1}(\mathbb{B}_\delta(0))$, where $\delta$ is chosen
 sufficiently small such that $\overline{U}_i$ does not contain any $x_j^\ast$ with $j\in\{1,\cdots,q\}\setminus\{i\}$.
 It follows from Propositions 5.1 and 5.8 that for some integer
 $I^{(i)}$ there holds
 $$\lim_{L\ra\infty}\lim_{k\ra\infty}\int_{\phi_i^{-1}(\mathbb{B}_{\delta\,/L}(0))}u_kf_k(x,u_k)dv_g=4\pi I^{(i)}.$$
 By Proposition 3.1, $u_k\ra u_\infty$ in $C^1_{\rm
 loc}(\Sigma\setminus\{x_1^\ast,\cdots,x_q^\ast\})$ as $k\ra\infty$.
 hence
 $$\lim_{L\ra\infty}\lim_{k\ra\infty}\int_{\Sigma\setminus \cup_{i=1}^q\phi_i^{-1}(\mathbb{B}_{\delta\,/L}(0))}
 u_kf_k(x,u_k)dv_g=\int_\Sigma u_\infty f_k(x,u_\infty)dv_g.$$
 Combining these two estimates, we obtain
 $$\lim_{k\ra\infty}\int_{\Sigma}
 u_kf_k(x,u_k)dv_g=\int_\Sigma u_\infty f_k(x,u_\infty)dv_g+4\pi\sum_{i=1}^qI^{(i)}.$$
 This together with (\ref{e1}) leads to
 $$\lim_{k\ra\infty}\int_{\Sigma}(|\nabla_gu_k|^2+\tau_k u_k^2)dv_g
 =\int_{\Sigma}(|\nabla_gu_\infty|^2+\tau_\infty u_\infty^2)dv_g
 +4\pi\sum_{i=1}^qI^{(i)}.$$
 In view of (\ref{0.6}), or particularly (\ref{F}), we then have
 $$\lim_{k\ra\infty}J_k(u_k)=J_\infty(u_\infty)+4\pi\sum_{i=1}^qI^{(i)}.$$
 This completes the proof of Theorem 1.1. $\hfill\Box$\\

 {\bf Acknowledgements.} This work is supported by the NSFC
 11171347. The author thanks the referee for his careful reading and
  valuable suggestions on the first version of this paper.

\end{document}